\documentclass{article}
\usepackage{multirow}
\usepackage{graphicx} % handles graphics and figures
\usepackage{mathtools,bm}
\usepackage{bbold}   %Get fancy double struck math notation for sets
\usepackage{amsmath}
\usepackage{courier}   %Have code written out nicely
\usepackage{color}
\usepackage{tabu}
\usepackage{soul}
\usepackage{algorithm}
\usepackage[noend]{algpseudocode}
\usepackage{float}
\usepackage[top=1in, bottom=1in, left=1in, right=1in]{geometry}
\usepackage[titletoc,title]{appendix}
\usepackage{amsthm}
\usepackage[colorlinks,
linkcolor=blue,
anchorcolor=blue,
citecolor=blue]{hyperref}
\usepackage{verbatim}
\numberwithin{equation}{section}
\newcommand{\bs}{\boldsymbol}

\definecolor{lightblue}{rgb}{0.0,0.5,0.8}

\definecolor{darkgreen}{rgb}{0.0, 0.4, 0.2}
\newtheorem{theorem}{Theorem}[section]
\newtheorem{lemma}{Lemma}[section]

\makeatother
\title{\textbf{Solutions of Nonlinear Optimal Control Problems Using Quasilinearization and Fenchel Duality}}
\author{Hailing Wang, Di Wu, Changjun Yu}

\begin{document}
\maketitle
\section*{Abstract}
\qquad In this paper, we consider a special class of nonlinear optimal control problems, where the control variables are box-constrained and the objective functional is strongly convex corresponding to control variables and separable with respect to the state variables and control variables.
We convert solving the original nonlinear problem into solving a sequence of constrained linear-quadratic optimal control problems by quasilinearization method. In order to solve each linear-quadratic problem efficiently we turn to study its dual problem. We formulate dual problem by the scheme of Fenchel duality, the strong duality property and the saddle point property corresponding to primal and dual problem are also proved, which together ensure that solving dual problem is effective. Thus solving the sequence of control constrained linear-quadratic optimal control problems obtained by quasilinearization technique is substituted by solving the sequence of their dual problem. We solve the sequence of dual problem and obtain the solution to primal control constrained linear-quadratic problem by the saddle point property. Furthermore, the fact that solution to each subproblem finally converges to the solution to the optimality conditions of original nonlinear problem is also proved. After that we carry out numerical experiments using present approach, for each subproblem we formulate the discretized primal and dual problem by Euler discretization scheme in our experiments. Efficiency of the present method is validated by numerical results.

\section{Introduction}
\qquad Optimal control is a subject that aims at controlling a given dynamic system over a period of time such that a specified performance index is minimized while any other constraints are satisfied in the process. Optimal control problems are widely encountered as mathematical models in such areas as industrial engineering \cite{MR1610543,MR2535357}, medical science \cite{MR4167193} and aerospace science \cite{MR3864034}. 

Nonlinear optimal control problems subjected to continuous or pointwise constraints, which are imposed at every moment along state and control variables have been an active area of research over the last several decades because of its wide applications. The main theoretical result for solving these problems analytically is minimum principle \cite{MR836572,MR0189864} which is regarded as first order necessary conditions. But the application of minimum principle to these problems finally results in a system of coupled two-point boundary value problem whose analytical solution is very difficult to obtain. Thus numerical methods are indispensible for solving applied nonlinear optimal control problems.

One of the most popular numerical methods is to approximate the problem by some discretization technique, such as using Euler discretization, precisely, using Euler scheme to discrete the dynamic constraints and objective functional, and finally obtain a large-scale finite dimensional optimization problem. Then solving original problem is converted to solving a finite dimensional optimization problem to get an approximate solution to original nonlinear infinity dimensional problem. Many finite dimensional optimization methods can be applied to solve the discretized problem such as SQP \cite{MR1713114}, interior point method \cite{MR1713114}, etc. This technique for original nonlinear optimal control problems is classified as direct method since it firstly discrete then optimize. Besides convergence theory for a wide range of optimal control problems solved by this class method is well established \cite{MR1779739, MR1786339, MR2208980}. 

The above numerical methods can be treated as solving the nonlinear problem directly. Futhermore we can also solve the original problem by adopting quasilinearization technique \cite{MR0178571,MR724931}. Then solving the nonlinear optimal control problem is converted into solving a sequence of linear-quadratic optimal control problems and each problem is much easier to solve. Then the critical problem is replaced by how to solve these linear-quadratic optimal control problems efficiently. Traditionally, each subproblem is solved by employing the aforementioned direct method such as Euler method \cite{MR1804658}, pseudospectral method \cite{MR2878201}, control parameterization method \cite{MR3119153}  and finally obtain a large-scale quadratic programming problem. We notice that  
Burachik and Kaya and Majeed have derived the dual problem of the control constrained linear-quadratic optimal control problem by Fenchel duality scheme in \cite{MR3199411}, the strong duality and saddle point property are also discussed in their paper, which guaranteed that solving dual problem is effective. Besides they illustrated some examples and showed via numerical experiments that by solving
the dual of the linear-quadratic problem, instead of the primal one, sometimes significantly improve the computational efficiency. We refer the reader to \cite{MR3199411} about the solving technique based on duality theory for control constrained linear-quadratic optimal control problems and the seminal papers of Rockafellar \cite{MR266020, MR282283, MR885199} for duality theory of linear-quadratic optimal control problems.

In this paper we develop a new iteration method for solving control-constrained nonlinear optimal control problems by combining quasilinearization technique with the duality theory for control-constrained linear-quadratic optimal control problems. In the proposed method, the nonlinear problem is converted into a series of control constrained linear-quadratic problem, then every subproblem's dual problem can be formulated, the solution to each subproblem is obtained by solving its dual problem using Euler scheme and taking advantage of the saddle point property. The method proposed is demonstrated efficiently on several practical examples and also compared to first quasilinearization then directly solve the sequence of subproblems. Numerical results show that present method behaves efficiently.

The rest of paper is organized as follows. In section 2, we make some preliminaries for our discussion. In section 3, we formulate the nonlinear control-constrained optimal control problem considered in this paper. We discuss the quasilinearization technique and convert solving the original nonlinear problem into solving a sequence of linear-quadratic problem. In section 4 and section 5, we follow the techniques used in refercence \cite{MR3199411} to form the dual problem of the general control-constrained linear-quadratic problem. Besides the linear-quadratic problem considered by us is more general than \cite{MR3199411}, the quadratic state term in the cost functional
is only positive semidefinite so that two cases considered in \cite{MR3199411} can be treated as special cases. Besides, the linear dynamic system is also more genernal camparing to \cite{MR3199411}. We derive dual problem and prove strong duality
and saddle point property. In section 6, we propose the new algorithm called sequential dual method. The convergence of algorithm is discussed in section 7. In section 8, we carry out numerical experiments with five examples to demonstrate present method is efficient. Finally, we make some conclusion and propose some outlook for our future work in section 9.

\section{Preliminaries}

\qquad Let $\mathbb{R}^n$ denotes the n-dimensional Euclidean space with the inner product given by $\bs x^{\top}\bs y$ and the norm
\begin{equation*}
\|\bs x\|=\sqrt{\bs x^{\top}\bs x}
\end{equation*}

Let $t_0,t_f\in \mathbb{R}$ and $t_0<t_f$. We denote by $L^2( t_0,t_f;\mathbb{R}^n)$ be the Hilbert space of Lebesgue measurable functions $\bs z:\lbrack t_0,t_f\rbrack\rightarrow \mathbb{R}^n$, with the inner product
\begin{equation*}
\langle \bs u,\bs v\rangle =\int_{t_0}^{t_f} \bs u^{\top}\bs v dt, \quad\text{for all\quad} \bs u,\bs v\in L^{2}(t_0,t_f;\mathbb{R}^n)
\end{equation*}
and equip it with $L^2$ norm, precisely,
\begin{equation*}
L^2( t_0,t_f;\mathbb{R}^n)=\left\{\bs z:\lbrack t_0,t_f\rbrack\rightarrow \mathbb{R}^n|\|\bs z\|_2=\left(\int_{t_0}^{t_f}|\bs z(t)|^2dt\right)^{1/2}<\infty\right\}
\end{equation*}
where $|\cdot|$ denotes the modulus of vector.

And $L^1( t_0,t_f;\mathbb{R}^n)$ is the Banach space of Lebesgue measurable functions $\bs z:\lbrack t_0,t_f\rbrack\rightarrow \mathbb{R}^n$
equipped with $L^1$ norm, precisely,
\begin{equation*}
L^1( t_0,t_f;\mathbb{R}^n)=\left\{\bs z:\lbrack t_0,t_f\rbrack\rightarrow \mathbb{R}^n|\|\bs z\|_1=\int_{t_0}^{t_f}|\bs z(t)|_1dt<\infty\right\}
\end{equation*}
where $|\cdot|_1$ denotes the one norm of vector.

$C( t_0,t_f;\mathbb{R}^n)$ is the Banach space of continuous vector functions $\bs z:\lbrack t_0,t_f\rbrack\rightarrow \mathbb{R}^n$ equipped with norm
\begin{equation*}
\|\bs z\|_{\infty}=\max\limits_{1\leq i\leq n}\{\mathop{\sup}\limits_{t\in\lbrack t_0,t_f\rbrack}|z_i(t)|\}
\end{equation*}

While $L^\infty( t_0,t_f;\mathbb{R}^n)$ is the Banach space of essentially bounded vector functions
with the norm
\begin{equation*}
\|\bs z\|_{\infty}=\max\limits_{1\leq i\leq n}\{\mathop{ess\sup}\limits_{t\in\lbrack t_0,t_f\rbrack}|z_i(t)|\}
\end{equation*}

Furthermore, $W^{1,2}(t_0,t_f;\mathbb{R}^n)$ is the Sobolev space of absolutely continuous functions, namely,
\begin{equation*}
W^{1,2}(t_0,t_f;\mathbb{R}^n)=\left\{\bs z\in L^2(t_0,t_f;\mathbb{R}^n)|\dot{\bs z}=\frac{d\bs z}{dt}\in L^2(t_0,t_f;\mathbb{R}^n)\right\}
\end{equation*}

\section{Nonlinear Optimal Control Problem with Control Constraints}
\qquad We consider the following nonlinear optimal control problem denoted by (NP):
\begin{equation*}
\text{(NP)}\left\{
\begin{aligned} 
\min\limits_{\bs x,\bs u}\quad &\int_{t_0}^{t_f} f(\bs x)+g(\bs u) dt \label{nonlinear.prob}\\
\text{subject to}\quad &\dot{\bs x}=\bs h(\bs x,\bs u)\quad\text{for }t \in \lbrack t_0, t_f\rbrack,\\
&\bs x(t_0)=\bs x_{0}\quad, E\bs x(t_f)=\bs e_f,\\
&\bs u\in U=\{\bs u\in L^2\lbrack t_0,t_f\rbrack|\alpha(t)\leq \bs u(t)\leq \beta(t)\quad\text{for all }t \in \lbrack t_0, t_f\rbrack\}
\end{aligned}
\right.
\end{equation*}
where the time horizon $\lbrack t_0,t_f\rbrack$ is specified, $\bs x\in W^{1,2}(t_0,t_f;\mathbb{R}^n),\bs u\in L^2(t_0,t_f;\mathbb{R})$ are state variables and control variable, respectively. And vector function $\bs h$ is assumed to be continuously differentiable with respect to $(\bs x,\bs u)$. It is assumed that the nonlinear dynamic system is controllable thus the problem is feasible. Note that the objective functional is separable corresponding to state variables $\bs x$ and control variable $\bs u$. Furthermore we assume function $f$ is convex and second-order continuously differentiable, function $g$ is strongly convex and second-order continuously differentiable with respect to their variables. Now the problem is to find the optimal control $\bs u$ and the corresponding state $\bs x$ satisfying the constraints while minimizing the objective functional. In this paper, we only consider a single control variable but our approach can be extended to the case with more than one control variable without much effort.
\subsection{Quasilinearization Technique for Problem(NP)}
\qquad Quasilinearization method is the generalization of Newton-Raphson method to functional space. Bellman and Kalaba \cite{MR0178571} pointed out that quasilinearization technique can be applied to solve variational problems in two different ways. First way is to linearize the differential equations derived by necessary conditions, namely, the optimality conditions, with boundary conditions. Thus the nonlinear variational problem can be solved by solving a sequence of linear differential equations with two-point boundary conditions. The second way is to expand the objective functional to second-order term and linearize the dynamic system and other nonlinear constraints around nominal state varables. In this way the nonlinear variational problem can be solved successively a sequence of linear-quadratic variational problems. Bashein and Enns \cite{quasi} applied quadratic programming to solve the nonlinear optimal control problems through quasilinearization technique. Jaddu \cite{MR1931507} took advantage of quasilinearization method to convert solving the nonlinear optimal control problem into solving a sequence of quadratic programming problems via state parameterization by Chebyshev series with unknown coefficients. Mohammad Maleki and Ishak Hashim \cite{MR3151790} combined  pseudospectral method with quasilinearization technique to solve the constrained time-delay nonlinear optimal control problems.  

Here we apply the idea of the second way of quasilinearization to the control constrained optimal control problem(NP). Expanding the objective functional and nonlinear dynamic system of problem(NP) around nominal state $\bs{x^N}$ and nominal control $\bs{u^N}$, the following control constrained linear-quadratic optimal control problem can be derived. We denote the subproblem as problem $Q^{N+1}$.
\begin{equation*}
\text{($Q^{N+1}$)}\left\{
\begin{aligned}
\min\limits_{\bs {x},\bs{u}}\quad &\int_{t_0}^{t_f} \frac{1}{2}(\bs {x}-\bs {x^{N}})^{\top}\nabla ^{2}f(\bs {x^{N}})(\bs {x}-\bs {x^{N}})+
\nabla f(\bs {x^{N}})(\bs {x}-\bs {x^{N}})+f(\bs {x^{N}})\\
&+\frac{1}{2}g''(\bs{ u^{N}})(\bs{u}-\bs {u^{N}})^{2}+ g'(\bs {u^{N}})(\bs{u}-\bs {u^{N}})+g(\bs {u^{N}})dt\\
\text{subject to}\quad &\bs{\dot{x}}=\bs h(\bs {x^{N}},\bs{ u^{N}})+
\bs h_{\bs x}(\bs{ x^{N}},\bs{ u^{N}})(\bs{x}-\bs{ x^{N}})+
\bs h_{\bs u}(\bs {x^{N}},\bs {u^{N}})(\bs{u}-\bs {u^{N}})\\
&\bs x(t_0)=\bs x_{0}\quad, E\bs x(t_f)=\bs e_f\\
&\bs {u}\in U=\{\bs u\in L^2\lbrack t_0,t_f\rbrack|\alpha(t)\leq \bs u(t)\leq \beta(t)\quad\text{for }t \in \lbrack t_0, t_f\rbrack\}
\end{aligned}
\right.
\end{equation*}

The procedure for solving the sequence of control constrained linear-quadratic optimal control problems starts with setting $N=0$ and choosing initial guess $\bs{x^0}$ and $\bs{u^0}$ which should not be too poor that cause quasilinearization method doesn't work. Then solving the obtained linear-quadratic optimal control problem $Q^1$ using numerical method such as pseudospectral method, Euler method, control parameterization method etc, denote the solution of problem $Q^1$ as $(\bs{\bar{x}^1},\bs{\bar{u}^1})$
%by converting it into dual problem as discussed in section 4.2, using the saddle point property given by Theorem 4.2, 
 Combining $(\bs{x^0},\bs{u^0})$ with $(\bs{\bar{x}^1},\bs{\bar{u}^1})$ a new nominal state and control can be constructed denoted as $(\bs{x^1},\bs{u^1})$. Then the next linear-quadratic optimal control problem $Q^2$ can be formulated by $(\bs{x^1},\bs{u^1})$. Repeat the procedure, thus, the original nonlinear problem is solved by solving a sequence of linear-quadratic optimal control problems. 

Thus efficiently solving each subproblem $Q^{N+1}$ obtained by quasilinearization technique plays important role on solving the nonlinear problem (NP) in our scheme. We note that Burachik and Kaya and Majeed studied the dual problem of some special linear-quadratic optimal control problem in \cite{MR3199411} and they show that in some case solving dual problem is more efficient than solving primal problem. Thus we consider the dual problem of the linear-quadratic problem of the form $Q^{N+1}$ in the next section. 

\section{The primal linear-quadratic optimal control problem}
\subsection{Problem Statement}
\qquad In this section, we consider the following class of linear-quadratic optimal control problem called problem (P), it is specfically that problem $(Q^{N+1})$ is a special case of (P). 

$$\text{(P)}\left\{
\begin{aligned}
 \min\limits_{\bs x,\bs u}\quad &\int_{t_0}^{t_f}\frac{1}{2}\bs x^{\top}\bs W(t)\bs x +\bs \omega(t)^{\top} \bs x+\frac{R(t)}{2}\cdot \bs u^{2}+r(t)\cdot \bs u dt\\
 \mbox{subject to}\quad &\dot{\bs x} =\bs A(t) \bs x+\bs B(t) \bs u +\bs c(t)\quad \text{for }t \in \lbrack t_0, t_f\rbrack,\\
  \quad &\bs x(t_0) = \bs x_0,  \quad  E\bs x(t_f) = \bs e_f, \\
  &\bs u(t)\in U(t)=\lbrack\alpha(t),\beta(t) \rbrack\quad \text{for }t \in \lbrack t_0, t_f\rbrack
\end{aligned}
\right.$$

The time horizon is taken to be $\lbrack t_0, t_f\rbrack$ with $t_0$ and $t_f$ specfied. The state variable $\bs x\in W^{1,2}(t_0,t_f;\mathbb{R}^n)$ and control variable $\bs u\in L^2(t_0,t_f;\mathbb{R})$. The time-varying matrices $\bs A : \lbrack t_0, t_f\rbrack\rightarrow \mathbb{R}^{n\times n}$ and $\bs B : \lbrack t_0, t_f\rbrack\rightarrow \mathbb{R}^{n}$ are continuous, $\bs W : \lbrack t_0, t_f\rbrack\rightarrow \mathbb{R}^{n\times n}$ is semi-positive definite and continuous, $\bs \omega,\bs c : \lbrack t_0, t_f\rbrack\rightarrow \mathbb{R}^{n}$ is continuous, and $R : \lbrack t_0, t_f\rbrack\rightarrow \mathbb{R}$ is positive and continuous, $r : \lbrack t_0, t_f\rbrack\rightarrow \mathbb{R}$ is continuous.

The initial state $\bs x_0\in \mathbb{R}^n$, matrix $E\in\mathbb{R}^{k\times n}$ and $\bs e_f\in \mathbb{R}^k$ are specified. 

The functions $\alpha,\beta : \lbrack t_0, t_f\rbrack\rightarrow \mathbb{R}$ are continuous.The constraint $\bs u(t)\in U(t)\quad (\text{for }t \in \lbrack t_0, t_f\rbrack)$ is sometimes called box-constraint (is to be satisfied pointwise).

The feasible set $\mathcal{F}$ is given by
\begin{align*}
\mathcal{F}=\{&(\bs x,\bs u)\in W^{1,2}(t_0,t_f;\mathbb{R}^n)\times L^2( t_0,t_f;\mathbb{R})|\\
&\dot{\bs x} =\bs A(t) \bs x+\bs B(t) \bs u +\bs c(t), \bs u(t)\in U(t)\quad \text{for }t \in \lbrack t_0, t_f\rbrack, \bs x(t_0) = \bs x_0,  \quad  E\bs x(t_f) = \bs e_f\}
\end{align*}
We will assume the feasible set $\mathcal{F}$ is nonempty. Note that $\mathcal{F}$ is convex and closed subset with respect to space $W^{1,2}(t_0,t_f;\mathbb{R}^m)\times L^2( t_0,t_f;\mathbb{R})$, and the objective functional is convex and continuous over $\mathcal{F}$, any local minimizer is also global minimizer.

We will assume that solution to problem(P) exists.

Remark 1.

(1) The control constraint set $U(t)$ can also be defined in terms of a one-side bound or no bound, i.e. , the following set may appear
\begin{equation*}
U(t)=\lbrack\alpha(t),+\infty)\quad,\quad U(t)=(-\infty,\beta(t)\rbrack\quad \text{or}\quad U(t)=(-\infty,+\infty)
\end{equation*}

(2) For simplicity in appearance, in the rest of this section, we omit in the calculations the argument $t$ of the functions $\bs W,\bs \omega,R,r,\bs A,\bs B,\bs c$, the variables $\bs x,\bs u$, and the bounds $\alpha,\beta$, whenever appropriate.
\subsection{Optimality Conditions}
\qquad In this section, we derive the optimality conditions for problem(P). We do this by means of the Pontryagin minimum principle. Since problem(P) is convex and continuous corresponding to variables $\bs x,\bs u$, the necessary conditions are also sufficient \cite{MR0189864}. Let $(\bs x,\bs u)\in\mathcal{F}$ be a solution to problem(P). We define the Hamiltonian function associated with problem(P) as follows: 
\begin{equation*}
H(\bs x,\bs u,\bs \lambda)=\lambda_0(\frac{1}{2}\bs x^{\top}\bs W\bs x +\bs\omega^{\top}\cdot \bs x+\frac{R}{2}\cdot \bs u^{2}+r\cdot \bs u)+\bs\lambda^{\top}(\bs A \bs x+\bs B \bs u +\bs c)
\end{equation*}
where the $\bs\lambda(t)\in W^{1,2}(t_0,t_f;\mathbb{R}^n)$ is called costate vector and $\lambda_0$ is a nonnegative constant. Note that if the dynamic system of problem (P) is controllable then we have $\lambda_0>0$ i.e. the problem (P) is normal. Motivated by this result, we assume that the dynamic system of problem (P) is controllable, and we can set constant $\lambda_0=1$ without loss of generality. 

The optimality conditions 
\begin{align}
\bs{\dot{\lambda}}&=-H_{\bs x}=-\bs W\bs x-\bs A^{\top}\bs\lambda-\bs\omega \label{adjoint.eq}\\
\bs \lambda(t_f)&=E^{\top}\bs z\label{adjoint.ter.}\\
\bs u&=\mathop{\arg\min}\limits_{\bs v\in U}H(\bs x,\bs v,\bs\lambda)
\end{align}
Here $\bs z\in \mathbb{R}^{k}$ denotes the multiplier corresponding to terminal state constraint $E\bs x(t_f) = \bs e_f$.
 
Since control variable is one-dimensional, the closed form of control variable can directly be specified.
\begin{equation}\label{optimal_control}
\bs u=\left\{
\begin{aligned}
&\frac{\bs B^{\top} (-\bs \lambda)-r}{R}\quad&\text{if}\quad\alpha\leq\frac{\bs B^{\top} (-\bs \lambda)-r}{R}\leq\beta,\\
&\alpha\quad&\text{if}\quad \frac{\bs B^{\top} (-\bs \lambda)-r}{R}<\alpha,\\
&\beta\quad&\text{if} \quad\frac{\bs B^{\top} (-\bs \lambda)-r}{R}>\beta.
\end{aligned}
\right.
\end{equation}

It follows from (\ref{optimal_control}) that the optimal control $\bs u$ can be treat as 
\begin{equation*}
\bs u=\textbf{Pr}_{[\alpha,\beta]}(\frac{\bs B^{\top} (-\bs \lambda)-r}{R})
\end{equation*}
where $\textbf{Pr}$ denotes the projection operator.  

Recall that $\alpha,\beta$ are continuous, therefore the optimal control $\bs u$ corresponding to problem(P) is continuous, so $\bs u$ belongs to $L^2(t_0,t_f;\mathbb{R})$. Substituting the optimal control $\bs u$(\ref{optimal_control}) in the state equations, we get the following differential equations:
\begin{equation}
{\bs{\dot x}} =\bs A \bs x+\bs B \bs u +\bs c
\end{equation}
with the boundary conditions
\begin{equation}\label{b.c.}
\bs x(t_0)=\bs x_0,\quad E\bs x(t_f)=\bs e_f
\end{equation}

Equations (\ref{optimal_control})-(\ref{b.c.}) together with (\ref{adjoint.eq})-(\ref{adjoint.ter.}) constitute the optimality conditions for problem(P)

Remark 2.

(1) If $E$ is inversable, then we note that terminal state $\bs x(t_f)$ can be specified. If the terminal state is unconstrained, we can set $E=0$ and $\bs e_f=0$ this case the costate vector satisfy $\bs\lambda(t_f)=0$. 

\section{The Dual Problem of LQ}
\subsection{Reformulation of Problem(P)}
\qquad We formulate the dual problem of problem(P) using Fenchel's duality scheme. In order to construct dual problem, we should firstly reformulate problem(P)
 such that the feasible set of new problem is a subspace. For this purpose, we can introduce some artificial variables $\bs s_1,\bs s_2\in\mathbb{R}^n$ and incorporate the boundary conditions and control constraint into the objective functional as indicator functions. Then we can rewrite the problem(P) equivalently as follows:
 \begin{align*}
\min\limits_{\bs x,\bs u,\bs s_1,\bs s_2}\quad&\int_{t_0}^{t_f}\frac{1}{2}\bs x^{\top}\bs W\bs x +\bs\omega^{\top} \bs x+\frac{R}{2}\cdot \bs u^{2}+r\cdot \bs u dt +\delta_{\bs x_0}(\bs s_1)+\delta_{\bs e_f}( E\bs s_2)+\delta_{1}(h)+\delta_{U}(\bs u)\\
\mbox{subject to}\quad&\dot{\bs x} =\bs A \bs x+\bs B \bs u +h\cdot\bs c \quad\text{for}\quad t \in \lbrack t_0, t_f\rbrack,\\ 
\quad &\bs x(t_0) = \bs s_1,  \quad  E\bs x(t_f) =E \bs s_2.
\end{align*}
where $\delta_{C}(a)$ is an indicator function of set $C$, which has value 0 when $a$ belongs to set $C$ and value $+\infty$ when $a$ doesn't belong to $C$. This case the feasible set $S$ is
\begin{equation}\label{feasible_set}
S=\{(\bs x,\bs u,h,\bs s_1,\bs s_2)|\dot{\bs x} =\bs A \bs x+\bs B\bs u +h\cdot\bs c\quad \text{for}\quad t \in \lbrack t_0, t_f\rbrack, \quad \bs x(t_0) = \bs s_1,  \quad  E\bs x(t_f) =E \bs s_2\}
\end{equation}
Note that feasible set $S$ is now a closed subspace of $W^{1,2}(t_0,t_f;\mathbb{R}^n)\times L^2( t_0,t_f;\mathbb{R})\times\mathbb{R}\times\mathbb{R}^n\times\mathbb{R}^n$. By the means of the definition of indicator function it is easy to show the equivalence between the rewritten problem and problem(P).

The objective functional of the rewritten problem is separable in the variables $\bs x,\bs u,\bs s_1,\bs s_2$. Define the following functionals:
\begin{align*}
f_1(\bs x) &=\int_{t_0}^{t_f}\frac{1}{2}\bs x^{\top}\bs W\bs x +\bs\omega^{\top} \bs xdt \\
f_2(\bs u)&=\int_{t_0}^{t_f}\frac{R}{2}\bs u^{2}+r\cdot\bs u dt+\delta_U(\bs u)
\end{align*}
as well as

$f_3(\bs s_1)=\delta_{\bs x_0}(\bs s_1)$,\quad $f_4(\bs s_2)=\delta_{\bs e_f}(E\bs s_2)$\quad and \quad$f_5(h)=\delta_{1}(h)$\\

Then objective functional can be treated as $f_1(\bs x)+f_2(\bs u)+f_3(\bs s_1)+f_4(\bs s_2)+f_5(h)$

\subsection{Derivative of the Dual Problem}
\qquad Following the Fenchel's duality scheme \cite{MR3204129}, we formulate the following dual problem of the problem(P). 
\begin{equation}\label{dual_P}
\begin{aligned}
\max\limits_{\bs x^*,\bs u^*,h^*,\bs s_1^*,\bs s_2^*}\quad &-(f_1^{*}(\bs {x^{*}})+f_2^{*}(\bs {u^{*}})+f_3^{*}(\bs {s_1^{*}})+f_4^{*}(\bs {s_2^{*}})+f_5^{*}(h^{*}))\\
\mbox{subject to} \quad &(\bs {x^{*}},\bs {u^{*}},h^{*},\bs {s_1^{*}},\bs{ s_2^{*}})\in S^{\perp}.
\end{aligned}
\end{equation}
where 
\begin{align*}
	f_{1}^{*}\left(\bs x^{*}\right) &=\sup _{\bs x \in W^{1,2}\left(t_{0}, t_{f} ; \mathbb{R}^{n}\right)}\left\{\langle\bs x^{*}, \bs x\rangle-f_{1}(\bs x)\right\} \\
	f_{2}^{*}\left(\bs u^{*}\right) &=\sup _{\bs u \in L^{2}\left(t_{0}, t_{f} ; \mathbb{R}\right)}\left\{\langle \bs u^{*},\bs u\rangle-f_{2}(\bs u)\right\} \\
	f_{3}^{*}\left(\bs s_{1}^{*}\right) &=\sup _{\bs s_{1} \in \mathbb{R}^{n}}\left\{\left\langle\bs s_{1}^{*}, \bs s_{1}\right\rangle-f_{3}\left(\bs s_{1}\right)\right\} \\
	f_{4}^{*}\left(\bs s_{2}^{*}\right) &=\sup _{\bs s_{2} \in \mathbb{R}^{n}}\left\{\left\langle\bs s_{2}^{*}, \bs s_{2}\right\rangle-f_{4}\left(\bs s_{2}\right)\right\}\\
	f_{5}^{*}\left(h^{*}\right) &=\sup _{h \in \mathbb{R}}\left\{\left\langle h^{*}, h\right\rangle-f_{5}\left(h\right)\right\}
\end{align*}
Futhermore $S^{\perp}\subseteq W^{1,2}(t_0,t_f;\mathbb{R}^n)\times L^2( t_0,t_f;\mathbb{R}^1)\times\mathbb{R}\times\mathbb{R}^n\times\mathbb{R}^n$ is the subspace orthogonal to $S$, i.e.
\begin{equation*}
S^{\perp}=\{(\bs x^*,\bs u^*,h^*,\bs s_1^*,\bs s_2^*)|\langle \bs x^*,\bs x\rangle+\langle \bs u^*,\bs u\rangle
+h^*h+\bs s_1^{*\top}\bs s_1+\bs s_2^{*\top}\bs s_2=0\}
\end{equation*}

To specify Fenchel dual problem, we should evaluate the conjugate functional $f_1^*,f_2^*,f_3^*,f_4^*,f_5^*$.

Note that the conjugate functional of $f_1$ is given by:
\begin{align}
f_1^{*}(\bs x^{*})&=\sup _{\bs x \in W^{1,2}\left(t_{0}, t_{f} ; \mathbb{R}^{n}\right)}\left\{\langle\bs x^{*}, \bs x\rangle-f_{1}(\bs x)\right\}\\
&=\sup\limits_{\bs x\in W^{1,2}\left(t_{0}, t_{f} ; \mathbb{R}^{n}\right)}\{\int_{t_0}^{t_f} (\bs x^{*}-\bs\omega)^{\top}\bs x-\frac{1}{2}\bs x^{\top}\bs W\bs x dt\}\label{f_1_con}
\end{align}
Now rewrite (\ref{f_1_con}) as
\begin{equation*}
-f_1^*(\bs x^{*})=\inf_{\bs x \in W^{1,2}\left(t_{0}, t_{f} ; \mathbb{R}^{n}\right)}\int_{t_0}^{t_f} F(\bs x,\bs x^*,t)dt
\end{equation*}
where $F(\bs x,\bs x^*,t)=-(\bs x^{*}-\bs\omega)^{\top}\bs x+\frac{1}{2}\bs x^{\top}\bs W\bs x$. Note that $\bs x$ solves the unconstrained minimization problem(\ref{f_1_con}) if and only if $\bs x$ solves the Euler-Lagrange equation \cite{MR3236927} 
\begin{equation*}
\frac{\partial F}{\partial x}-\frac{d}{dt}\left(\frac{\partial F}{\partial \dot{x}}\right)=0
\end{equation*}
Then it turns to be the equation where $\bs x$ is the solution to 
$\bs W\cdot \bs x+\bs\omega= \bs x^{*}$\\
Hence, we can obtain the closed form of the conjugate functional of $f_1$
\begin{equation}
f_1^{*}(\bs x^{*}) =\left\{
\begin{aligned}
\int_{t_0}^{t_f}\frac{1}{2}\bs y^{\top}\bs W\bs y dt\quad 
&\text{if}\quad  \bs W\cdot \bs y+\bs\omega= \bs x^{*}\\
+\infty\quad &\text{other cases}  
\end{aligned}
\right.
\end{equation}

Next, we compute the conjugate functional of $f_2$, defined by 
\begin{align}
	f_{2}^{*}\left(\bs u^{*}\right) &=\sup _{\bs u \in L^{2}\left(t_{0}, t_{f} ; \mathbb{R}\right)}\left\{\langle \bs u^{*},\bs u\rangle-f_{2}(\bs u)\right\} \\
	&=\sup\limits_{\bs u\in L^{2}\left(t_{0}, t_{f} ; \mathbb{R}\right)}\{\int_{t_0}^{t_f}(\bs u^*-r)\cdot\bs u-\frac{R}{2}\bs u^{2} dt-\delta_U(\bs u)\}\\
	&=\sup\limits_{\bs u\in U}\{\int_{t_0}^{t_f}(\bs u^*-r)\cdot\bs u-\frac{R}{2}\bs u^{2} dt\}\label{f_2_con}
\end{align}

This maximization problem can be treated as a linear-quadratic control problem with box constraints on the control variable. Since the objective functional is strongly concave, there exists a unique solution. We select $\bs u$ to minimize the integral term pointwise in $\lbrack t_0,t_f\rbrack$, namely, for every given $t\in\lbrack t_0,t_f\rbrack$,
\begin{equation*}
\bs u(t)=\mathop{\arg\max}\limits_{\alpha(t)\leq\bs v\leq\beta(t)}\{(\bs u(t)^*-r(t))\cdot\bs v-\frac{R(t)}{2}\bs v^{2}\}
\end{equation*}
Then solution is obtained as 
\begin{equation*}
\bs u=\textbf{Pr}_{[\alpha,\beta]}(\frac{\bs u^*-r}{R})
\end{equation*}
Substituting $\bs u$ into the (\ref{f_2_con}), we specify $f_2^*$
\begin{align}
f_2^{*}(\bs u^{*})=\int_{t_0}^{t_f}\psi (\bs u^{*}) dt
\end{align}
where $\psi (\bs u^{*})$ is
\begin{equation}\label{psi_def}
\psi(\bs u^{*})=\left\{
\begin{aligned}
&\frac{(\bs u^{*}-r)^{2}}{2\cdot R} \quad&\text{if}\quad\alpha\leq \frac{\bs u^{*}-r}{R}\leq \beta\\
&-\frac{R}{2}\cdot\alpha^{2}+(\bs u^{*}-r)\cdot \alpha
\quad&\text{if}\quad\frac{\bs u^{*}-r}{R}< \alpha\\
&-\frac{R}{2}\cdot\beta^{2}+(\bs u^{*}-r)\cdot \beta
\quad &\text{if}\quad\frac{\bs u^{*}-r}{R} > \beta
\end{aligned}
\right .
\end{equation}

The conjugate functional of $f_3$ is given by
\begin{equation}
f_3^{*}(\bs s_1^{*})=\sup\limits_{\bs s_1^*} \{\bs s_1^{\top}\bs s_1^*-\delta_{\bs x_0}(\bs s_1)\}
\end{equation}
which implies $\bs s_1=\bs x_0$, then $f_3^*(\bs s_1^*)=\bs x_0^{\top}\bs s_1^{*}$

The caculation of $f_4^*$ and $f_5^*$ is similar, finally we find $f_4^*$ and $f_5^*$ defined by
\begin{align}
&f_4^{*}(\bs s_2^{*})=\left\{
\begin{aligned}
\bs e_f^{\top}\bs \eta\quad &\text{if} \quad E^{\top}\bs\eta=\bs s_2^{*}\\
+\infty \quad &\text{other cases}
\end{aligned}
\right.\\
&f_5^{*}(h^{*})=h^{*}
\end{align}

Next we present, in Theorem 5.1, the orthogonal subspace of $S$, where $S$ is defined in (\ref{feasible_set}).  
\begin{theorem}\label{ortho_s}
     Consider the closed subspace $S$ of $W^{1,2}(t_0,t_f;\mathbb{R}^n)\times L^2(t_0,t_f;\mathbb{R})\times\mathbb{R}\times\mathbb{R}^n\times\mathbb{R}^n$ as described in $(\ref{feasible_set})$, then its orthogonal subspace $S^{\perp}$, is specified by
     \begin{align}
     S^{\perp}=\{ (\bs x^{*},\bs u^{*},h^{*},\bs s_1^{*},\bs s_2^{*})|\dot{\bs p}=-\bs A^{\top}\bs p+\bs x^{*},\bs p(t_f)=-\bs s_2^{*},\bs u^{*}=\bs B\bs p,\\
     \bs s_1^{*}=\bs p(t_0),h^{*}=\int_{t_0}
     ^{t_f}\bs p^{\top}\bs cdr \}\notag
     \end{align}
\end{theorem}
\emph{Proof}.
Since the proof is similar to \cite{MR3199411}(LEMMA 1), we omit proof here.
$\hfill\qedsymbol$

By Theorem \ref{feasible_set}, dual problem given by (\ref{dual_P}) can now be expressed explicitly as follows called problem (DP):
\begin{equation*}
\text{(DP)}\left\{
\begin{aligned} 
-\min \quad &\bs x_0^{\top}\bs p(t_0)+\bs e_f^{\top}\bs\eta+\int_{t_0}^{t_f}\frac{1}{2}\bs y^{\top}\bs W\bs y+\bs p^{\top}\bs c+\psi(\bs B^{\top}\bs p)dt\label{eq:dual}\\
\mbox{subject to} \quad  &\dot{\bs p}=-\bs A^{\top}\bs p+\bs W\bs y+\bs\omega \notag\\
&-\bs p(t_f)=E^{\top}\bs\eta\notag
\end{aligned}
\right.
\end{equation*}
The dual problem is a linear-quadratic optimal control problem, the state variable of problem(DP) is $\bs p\in W^{1,2}(t_0,t_f;\mathbb{R}^n)$ and control variable $\bs y\in L^2(t_0,t_f;\mathbb{R}^n)$. Comparing to problem(P), the control variable in problem (DP) is unconstrained i.e. $\bs y\in L^2(t_0,t_f;\mathbb{R}^n) $\\

Remark 3.

(1) The function $\psi$ defined in (DP) is continuous differentiable but not twice differentiable.

(2) Unlike the one-dimensional control variable $\bs u$ in the problem(P), the dual problem(DP)'s control variable $\bs y$ is n-dimensional. Precisely the dual problem (DP) compared with problem (P), it has more variables while less constraints. 

(3) Note that if $\bs W$ is positive-definite then the dynamic constraints in problem(DP) can be equivalently transformed into $\bs y=\bs W^{-1}(\dot{\bs p}+\bs A^{\top}\bs p-\bs\omega)$. Furthermore we can substitute it into objective functional which means that the problem can be equivalently transformed into a new problem without dynamic constraints.
\subsection{Optimality Conditions for Problem(DP) and Strong Duality}
\qquad In this section, we derive the necessary conditions for problem(DP). Since (DP) is continuous and convex in state variables $\bs p$ control variables $\bs y$. The optimality conditions are also sufficient. Since strong duality and saddle point property are important in any duality scheme, we prove the strong duality property for the problem(DP) and problem(P) as well as the saddle point property. 

Problem(DP) is also assumed to be normal, as discussed in section 4.2. Then the Hamiltonian function associated with problem(DP) is defined as follows:
\begin{equation*}
H(\bs p,\bs y,\bs \mu)=\frac{1}{2}\bs y^{\top}\bs W\bs y+\bs p^{\top}\bs c+\psi(\bs B^{\top}\bs p)+\bs\mu^{\top}(-\bs A^{\top}\bs p+\bs W\bs y+\bs\omega)
\end{equation*}
where $\bs\mu(t)\in W^{1,2}(t_0,t_f;\mathbb{R}^n)$ is called costate vector and function $\psi$ is defined by (\ref{psi_def}). Since the control variables $\bs y$ is unconstrained, then minimum principle $\bs y=\mathop{\arg\min}\limits_{\bs\gamma\in L(t_0,t_f;\mathbb{R}^n)}H(\bs p,\bs\gamma,\bs \mu)$ can be converted into
\begin{equation}\label{dp_condtion_1}
\bs W\bs y=-\bs W\bs \mu
\end{equation}
On the other hand, the costate equation for problem(DP) is given by
\begin{equation}
\bs{\dot\mu}=\bs A\bs\mu-\bs B\nabla_{\bs B^{\top}\bs p}\psi(\bs B^{\top}\bs p)-\bs c
\end{equation}
with the boundary conditions
\begin{equation}
\bs\mu(t_0)=-\bs x_0,\qquad E\bs\mu(t_f)=-\bs e_f
\end{equation}
Besides the state equation
\begin{equation}
\bs{\dot p}=-\bs A^{\top}\bs p+\bs W\bs y+\bs\omega
\end{equation}
with the boundary conditions
\begin{equation}\label{dp_condition_end}
-\bs p(t_f)=E^{\top}\bs \eta
\end{equation}

Equations (\ref{dp_condtion_1})-(\ref{dp_condition_end}) consititute the optimality conditions for the dual problem(DP).

\begin{theorem}
	Suppose that $(\bs x,\bs u)$ is the optimal solution to problem(P), and let the costate vector $\bs \lambda$, multiplier $\bs z$ be the solution to optimality conditions (\ref{adjoint.eq})-(\ref{adjoint.ter.}) corresponding to Problem(P).    Let 
	\begin{equation*}
	\quad\bs p=-\bs\lambda, \quad\bs \mu=-\bs x,\quad \bs y=\bs x, \text{and}\quad \bs \eta=\bs z 
	\end{equation*}
	Then there holds 
	
	(a) (Strong Duality) the optimal value of problem(P) and problem(DP) are equal.
	
	(b) (Saddle Point) $(\bs p,\bs y,\bs \mu,\bs \eta)$ defined by the solution to optimality conditions corresponding to the problem(P) also is a solution to optimality conditions corresponding to the problem (DP).
\end{theorem}
{\emph{Proof}}.
To prove the theorem, we firstly show that $(\bs p,\bs y)$ is feasible for problem(DP). By substituting $\bs p=-\bs \lambda$, $\bs \mu=-\bs x$, $\bs y=\bs x$ and $\bs \eta=\bs z$ into the conditions (\ref{adjoint.eq}) and (\ref{adjoint.ter.}), we get the equality constraints and the terminal state constraint in problem(DP), namely,
\begin{align*}
 \dot{\bs p}&=-\bs A^{\top}\bs p+\bs W\bs y+\bs\omega \\
 -\bs p(t_f)&=E^{\top}\bs\eta
\end{align*}
Hence $(\bs p,\bs y)$ is feasible for problem(DP). Next we should show the duality gap is zero, then the strong duality is obtained. Let $\pi(\bs x,\bs u)$, $\delta(\bs p,\bs y,\bs \eta)$ denote the objective functional value of primal problem(P) and dual problem(DP) respectively. Taking advantage of the definition of $(\bs p,\bs y,\bs\eta)$, we have
\begin{align*}
-\delta(\bs p,\bs y,\bs \eta)&=\bs x_0^{\top}\bs p(t_0)+\bs e_f^{\top}\bs\eta+\int_{t_0}^{t_f}\frac{1}{2}\bs y^{\top}\bs W\bs y+\bs p^{\top}\bs c+\psi(\bs B^{\top}\bs p)dt\\
&=-\bs x_0^{\top}\bs \lambda(t_0)+\bs x(t_f)^{\top}\bs\lambda(t_f)+\int_{t_0}^{t_f}\frac{1}{2}\bs x^{\top}\bs W\bs x-\bs \lambda^{\top}\bs c+\psi(\bs B^{\top}\bs p)dt\\
&=\int_{t_0}^{t_f}\frac{d}{dt}(\bs x^{\top}\bs \lambda) dt+\int_{t_0}^{t_f}\frac{1}{2}\bs x^{\top}\bs W\bs x-\bs \lambda^{\top}\bs c+\psi(-\bs B^{\top}\bs\lambda)dt\\
&=\int_{t_0}^{t_f}\frac{1}{2}\bs x^{\top}\bs W\bs x-\bs \lambda^{\top}\bs c+\psi(-\bs B^{\top}\bs\lambda)+\bs{\dot{x}}^{\top}\bs \lambda+\bs x^{\top}\bs {\dot{\lambda}} dt
\end{align*}
Furthermore $\bs\lambda$ satisfies the following equation 
\begin{equation*}
\bs{\dot{\lambda}}=-\bs W\bs x-\bs A^{\top}\bs\lambda-\bs\omega
\end{equation*}
Substituting this expression and rearranging, we get
\begin{equation}\label{dual_gap}
-\delta(\bs p,\bs y,\bs \eta)=\int_{t_0}^{t_f}-\frac{1}{2}\bs x^{\top}\bs W\bs x-\bs \omega^{\top}\bs x+\psi(-\bs B^{\top}\bs\lambda)+\bs \lambda^{\top}\bs B\bs u dt
\end{equation}
Since $\bs u$ is the optimal solution to problem(P), from (\ref{optimal_control}) and $\bs p=-\bs\lambda$, we have that
\begin{equation*}
\bs\lambda^{\top}\bs B\bs u=(\bs B^{\top}\bs\lambda) \textbf{Pr}_{\lbrack\alpha,\beta\rbrack}(-\frac{\bs B^{\top}\bs\lambda+r}{R})
\end{equation*}
From (\ref{psi_def}) we have that
\begin{equation*}
\psi(-\bs B^{\top}\bs\lambda)=-\frac{R}{2}(\textbf{Pr}_{\lbrack\alpha,\beta\rbrack}(-\frac{\bs B^{\top}\bs\lambda+r}{R}))^2-(\bs B^{\top}\bs\lambda+r)\textbf{Pr}_{\lbrack\alpha,\beta\rbrack}(-\frac{\bs B^{\top}\bs\lambda+r}{R})
\end{equation*}
Then we obtain 
\begin{align}
\psi(-\bs B^{\top}\bs\lambda)+\bs\lambda^{\top}\bs B\bs u&=-\frac{R}{2}(\textbf{Pr}_{\lbrack\alpha,\beta\rbrack}(-\frac{\bs B^{\top}\bs\lambda+r}{R}))^2-r\textbf{Pr}_{\lbrack\alpha,\beta\rbrack}(-\frac{\bs B^{\top}\bs\lambda+r}{R})\notag\\
&=-\frac{R}{2}\bs u^2-r\bs u\label{sum_function}
\end{align}
Substituting (\ref{sum_function}) into (\ref{dual_gap}) yields
\begin{equation*}
-\delta(\bs p,\bs y,\bs \eta)=\int_{t_0}^{t_f}-\frac{1}{2}\bs x^{\top}\bs W\bs x-\bs \omega^{\top}\bs x-\frac{R}{2}\bs u^2-r\bs u dt=-\pi(\bs x,\bs u)
\end{equation*}

Since we have shown $(\bs p,\bs y,\bs\eta)$ is feasible, the strong duality holds clearly. We must have that $(\bs p,\bs y,\bs\eta)$ is the optimal solution to problem(DP). Thus the proof is complete.
$\hfill\qedsymbol$
\section{Algorithm Design for Nonlinear Problem (NP)}
\qquad In this section, we design computational algorithm for solving problem (NP) that combines quasilinearization technique with the duality for linear-quadratic optimal control problems discussed in the previous section.

Firstly, we give algorithm for the case that terminal state is unconstrained. It can be described by problem (NP) by setting $E=0$ and $\bs e_f=0$ .
\begin{algorithm}[htb]
	\caption{Sequential Dual Method 1} %算法的名字
	\hspace*{0.02in} {\bf Step 1:} %算法的输入， \hspace*{0.02in}用来控制位置，同时利用 \\ 进行换行
	Input initial state $\bs{ x^0}$, initial control $\bs{u^0}$, stopping tolerence constant $tol$, iteration number $k=0$\\
	\hspace*{0.02in} {\bf Step 2:}
	Taking advantage of $\bs{ x^k}$ and $\bs {u^k}$ to formulate the quasilinearization subproblem $Q^{k+1}$\\
	\hspace*{0.02in} {\bf Step 3:} %算法的结果输出
	Formulate the dual problem of $Q^{k+1}$ and get the solution $(\bs {x^{k+1}},\bs {u^{k+1}})$ to problem $Q^{k+1}$ by the solution to its dual problem and saddle point property shown in Theorem 5.2.\\
	\hspace*{0.02in} {\bf Step 4:} %算法的结果输出
	Set $\bs{d^k}=(\bs {x^{k+1}}-\bs {x^{k}},\bs {u^{k+1}}-\bs {u^{k}})$.
	If $\|\bs{d^k}\|=\|\bs {x^{k+1}}-\bs {x^{k}}\|+\|\bs {u^{k+1}}-\bs {u^{k}}\|<tol$, then stop iteration and return $(\bs {x^{k+1}},\bs {u^{k+1}})$ as solution to original nonlinear problem. 
	
	Otherwise set $k=k+1$ and return \textbf{Step 2}.
\end{algorithm} 

For general case that the terminal state is constrained, we design the following algorithm.
\begin{algorithm}[htb]
	\caption{Sequential Dual Method 2} %算法的名字
	\hspace*{0.02in} {\bf Step 1:} %算法的输入， \hspace*{0.02in}用来控制位置，同时利用 \\ 进行换行
	Input initial state $\bs{ x^0}$, initial control $\bs{u^0}$, stopping tolerence constant $tol$, iteration number $k=0$\\
	\hspace*{0.02in} {\bf Step 2:}
	Taking advantage of $\bs{ x^k}$ and $\bs {u^k}$ to formulate the quasilinearization subproblem $Q^{k+1}$\\
	\hspace*{0.02in} {\bf Step 3:} %算法的结果输出
	Formulate the dual problem of $Q^{k+1}$ and get the solution $(\bs {\bar{x}^{k+1}},\bs {\bar{u}^{k+1}})$ to problem $Q^{k+1}$ by the solution to its dual problem and saddle point property shown in Theorem 5.2.\\
	\hspace*{0.02in} {\bf Step 4:} %算法的结果输出
	Set $\bs{d^k}=(\bs {\bar{x}^{k+1}}-\bs {x^{k}},\bs {\bar{u}^{k+1}}-\bs {u^{k}})$.
	If $\|\bs{d^k}\|=\|\bs{\bar{x}^{k+1}}-\bs {x^{k}}\|+\|\bs {\bar{u}^{k+1}}-\bs {u^{k}}\|<tol$, then stop iteration and return $(\bs {\bar{x}^{k+1}},\bs {\bar{u}^{k+1}})$ as solution to original nonlinear problem. 
	
	\hspace*{0.02in} {\bf Step 5:} %算法的结果输出
	By merit functional $P(\bs x,\bs u)$, select step length $\kappa_k$ such that  
	\begin{equation*}
	P((\bs{ x^k},\bs{u^k})+\kappa_k\bs{d^k})=\min\limits_{0\leq\kappa\leq 1}P((\bs{ x^k},\bs{u^k})+\kappa\bs{d^k})
	\end{equation*}
	Renew $(\bs{ x^{k+1}},\bs{u^{k+1}})=(\bs{ x^k},\bs{u^k})+\kappa_k\bs{d^k}$, $k=k+1$ then return \textbf{Step 2}.
\end{algorithm} 

Here we define merit functional $P(\bs x,\bs u)$ as following equation and $\theta$ is a constant.

\begin{equation*}
P(\bs x,\bs u)=\int_{t_0}^{t_f} f(\bs x)+g(\bs u) dt+\theta\|\dot{\bs x}-\bs h(\bs x,\bs u)\|_1
\end{equation*}
Note that the merit functional takes the original problem (NP) feasibility and the objective functional into account.

Comparing Algorithm 1 to Algorithm 2, we observe that Algorithm 2 has line-search procedure while Algorithm 1 doesn't have. We will analyze this in the next section in detail. 
\section{Convergence Analysis}
\qquad In this section, we prove the convergence of the sequence $\{(\bs{x^N},\bs{u^N})\}$ generated by sequential dual method designed in previous section for problem(NP).    

We assume that $\bs h$ is continuous corresponding to each argument, the dynamic system in original nonlinear problem (NP) is controllable. Furthermore we assume the dynamic system of every subproblem is controllable. 

Firstly, we shall prove that subproblem $Q^{N+1}$ obtained by quasilinearization using $(\bs{x^N},\bs{u^N})$ has at least one solution thus quasilinearization technique is reasonable for problem(NP).

\begin{lemma}
	Suppose that $\bs{x^N},\bs{u^N},\alpha,\beta$ are bounded, there exists some $\bs u\in L^2(t_0,t_f;\mathbb{R})$ be the optimal control for subproblem $Q^{N+1}$ 
\end{lemma}
{\emph{Proof}}. Let $J(\bs v)$ represents the objective functional, precisely, 
\begin{align*}
J(\bs v)=&\int_{t_0}^{t_f} \frac{1}{2}(\bs {y}-\bs {x^{N}})^{\top}\nabla ^{2}f(\bs {x^{N}})(\bs {y}-\bs {x^{N}})+
 f(\bs {x^{N}})^{\top}(\bs {y}-\bs {x^{N}})\\
&+\frac{1}{2}g''(\bs{ u^{N}})(\bs{v}-\bs {u^{N}})^{2}+ g'(\bs {u^{N}})(\bs{v}-\bs {u^{N}})dt
\end{align*}

Firstly we prove the objective functional of problem $Q^{N+1}$ is weakly lower semi-continuous in functional space $L^2(t_0,t_f;\mathbb{R})$.

Let $\{\bs{v_n}\}$ converges weakly to $\bs{\bar{v}}$ in $L^2(t_0,t_f;\mathbb{R})$ and $\bs {y_n}:=\bs x(t|\bs{v_n})$ represents the solution to the following dynamic system
\begin{align*}
&\bs{\dot{x}}=\bs h(\bs {x^{N}},\bs{ u^{N}})+
\bs h_{\bs x}(\bs{ x^{N}},\bs{ u^{N}})(\bs{x}-\bs{ x^{N}})+
\bs h_{\bs u}(\bs {x^{N}},\bs {u^{N}})(\bs{v_n}-\bs {u^{N}})\\
&\bs x(t_0)=\bs x_{0}, E\bs x(t_f)=\bs e_f
\end{align*}

Since $\bs {y_n}:=\bs x(t|\bs{v_n})$ is solution to linear differential equations, combining with the definition of weak convergence, we can directly derive that $\bs{y_n}$ converges to $\bs{\bar{y}}$ pointwisely, $\bs{\bar{y}}$ is the solution to the following differential equations,  
\begin{align*}
&\bs{\dot{x}}=\bs h(\bs {x^{N}},\bs{ u^{N}})+
\bs h_{\bs x}(\bs{ x^{N}},\bs{ u^{N}})(\bs{x}-\bs{ x^{N}})+
\bs h_{\bs u}(\bs {x^{N}},\bs {u^{N}})(\bs{\bar{v}}-\bs {u^{N}})\\
&\bs x(t_0)=\bs x_{0}, E\bs x(t_f)=\bs e_f
\end{align*}

Since function $g$ is strongly convex, $J(\bs v)$ can be represented as
\begin{align*}
J(\bs v)=&\int_{t_0}^{t_f} \frac{1}{2}(\bs {y}-\bs {x^{N}})^{\top}\nabla ^{2}f(\bs {x^{N}})(\bs {y}-\bs {x^{N}})+
\nabla f(\bs {x^{N}})^{\top}(\bs {y}-\bs {x^{N}})-\frac{1}{2}\frac{(g'(\bs{u^N}))^2}{g''(\bs{u^N})}dt\\
&+\frac{1}{2}\left\|\sqrt{g''(\bs{ u^{N}})}(\bs{v}-\bs {u^{N}})+\frac{g'(\bs {u^{N}})}{\sqrt{g''(\bs{ u^{N}})}}\right\|_2^2\\
\end{align*}

Because any norm of Banach space is weakly lower semi-continuous, taking advantage of Lebesgue-Dominated-Convergence theorem, we have that
$\lim\limits_{n\to\infty}\inf J(\bs{v_n})\geq J(\bs{\bar{v}})$ which means the objective functional $J$ is weakly semi-continuous.

Then, we can prove the existence of solution to problem $Q^{N+1}$. We observe that our objective functional $J$ is convex quadratic functional, the admissible control set is bounded thus the infimum of $J(\bs v)$ exists.

Assume that there is a minimizing sequence $\{\bs{v_n}\}$ such that 
\begin{equation*}
\lim\limits_{n\to\infty} J(\bs{v_n})=\inf\limits_{\bs v\in\lbrack \alpha(t),\beta(t)\rbrack}J(\bs v)
\end{equation*}
Since $\{\bs{v_n}\}$ is bounded, implies that there exists a subsequence, for simplicity, still denotes by $\{\bs{v_n}\}$, that converges weakly to some $\bs u$ in $L^2(\lbrack t_0,t_f\rbrack;\mathbb{R})$. It follows from weakly lower semi-continuous of $J$ and we thus have 
\begin{equation*}
J(\bs u)\leq\lim\limits_{n\to\infty}\inf J(\bs{v_n})=\inf\limits_{\bs v\in\lbrack \alpha(t),\beta(t)\rbrack}J(\bs v)
\end{equation*}
Hence, $\bs u$ must be the optimal control, which means that subproblem $Q^{N+1}$ exists at least one solution.
$\hfil\qedsymbol$\\

Since we have assumed the controllability condition for dynamic system of problem (NP), we can directly derive the optimality conditions for it. Suppose that $(\bs {x^{*}},\bs{ u^{*}})$ is the local minimizer of problem (NP).

Let $\bs{\lambda^*}$ denotes the solution to the linear costate equations and $\bs \eta^*$ denotes the multiplier corresponding to the terminal state constraints.
\begin{equation}\label{optimal_condtion_NP}
\bs{{\dot{\lambda^*}}}=-\nabla f(\bs{x^*})-\bs h_{\bs x}(\bs {x^{*}},\bs {u^{*}})^{\top}\bs{\lambda^{*}}
\end{equation}
with boundary conditions
\begin{equation}
\bs{ \lambda^{*}}(t_f)=E^{\top}\bs{\eta^*}
\end{equation}
The state equations
\begin{equation}
\bs{{\dot{x}}^{*}}=\bs h(\bs{ x^{*}},\bs {u^{*}})
\end{equation}
with boundary conditions
\begin{equation}
\bs {x^{*}}(t_0)=\bs x_{0}, E\bs {x^{*}}(t_f)=\bs e_f
\end{equation}

The Hamiltonian function with respect to this nonlinear problem is defined as
\begin{equation*}
H(\bs {x},\bs u,\bs{ \lambda})=f(\bs{x})+g(\bs u)+\bs{ \lambda}^{\top}h(\bs{x},\bs u)
\end{equation*}

From the minimum principle, it follows that
\begin{equation}\label{optimal_condition_np_end}
-H_{\bs u}(\bs {x^{*}},\bs {u^{*}},\bs{ \lambda^{*}})=-(g'(\bs{u^*})+\bs h_{\bs u}(\bs {x^{*}},\bs {u^{*}})^{\top}\bs{\lambda^{*}})\in \mathcal{N}_{U}(\bs {u^{*}})\quad\text{for all}\quad t\in\lbrack t_0,t_f \rbrack
\end{equation}
where $\mathcal{N}_{U}$ is the normal cone. Precisely, for any $t$ given in time horizon the normal cone of $\bs u$ at time point $t$ is defined by 
\begin{equation*}
\mathcal{N}_{U}(\bs u(t))=\{\omega\in\mathbb{R}|\omega\cdot(\bs v-\bs u(t))\leq 0, \text{for all}\quad \bs v\in {U(t)}\}
\end{equation*}
Then (\ref{optimal_condtion_NP})-(\ref{optimal_condition_np_end}) constitute the first-order optimality conditions for problem (NP).

As discussed in section 4.2, for problem $Q^{N+1}$ obtained by $(\bs{x^N},\bs{u^N})$ we can specify its optimality conditions.

Here we define the Hamiltonian function with respect to problem $Q^{N+1}$ as follows
\begin{align*}
H(\bs{ x},\bs{ u},\bs {\lambda})&=\frac{1}{2}(\bs x-\bs {x^{N}})^{\top}\nabla ^{2}f(\bs {x^{N}})(\bs x-\bs {x^{N}})+
\nabla f(\bs {x^{N}})^{\top}(\bs x-\bs {x^{N}})\\
&+\frac{1}{2}g''(\bs{ u^{N}})(\bs u-\bs {u^{N}})^{2}+ g'(\bs {u^{N}})(\bs u-\bs{ u^{N}})+\bs \lambda^{\top}(\bs h(\bs {x^{N}},\bs {u^{N}})\\
&+
\bs h_{\bs x}(\bs {x^{N}},\bs {u^{N}})(\bs x-\bs {x^{N}})+
\bs h_{\bs u}(\bs {x^{N}},\bs {u^{N}})(\bs u-\bs {u^{N}}))
\end{align*}

Let $\bs{\lambda^{N+1}}$ denotes the costate vector, the costate equation
\begin{equation}\label{QLNP_optimal_condtion_1}
\bs{\dot{\lambda}^{N+1}}=-\nabla^2f(\bs{ x^N})(\bs {x^{N+1}}-\bs {x^N})-\nabla f(\bs{ x^N})-\bs h_{\bs x}(\bs {x^{N}},\bs{ u^{N}})^{\top}\bs{\lambda^{N+1}}
\end{equation}
with boundary conditions
\begin{equation}\label{QLNP_con_cos}
\bs {\lambda^{N+1}}(t_f)=E^{\top}\bs{\eta^{N+1}}
\end{equation}
here $\bs{\eta^{N+1}}$ denotes the multiplier corresponding to the terminal state constraints.

The state equation
\begin{equation}\label{QLNP_dynamic}
\bs{\dot{ x}^{N+1}}=\bs h(\bs{ x^{N}},\bs {u^{N}})+
\bs h_{\bs x}(\bs{ x^{N}},\bs{ u^{N}})(\bs {x^{N+1}}-\bs {x^{N}})+
\bs h_{\bs u}(\bs {x^{N}},\bs{ u^{N}})(\bs {u^{N+1}}-\bs {u^{N}})
\end{equation}
with boundary conditions
\begin{equation}\label{dyna_con}
\bs {x^{N+1}}(t_0)=\bs x_{0}, E\bs {x^{N+1}}(t_f)=\bs e_{f}
\end{equation}
The first-order optimality conditions
\begin{equation}\label{QLNP_optimality_condition_end}
\bs -H_{\bs u}(\bs{ x^{N+1}},\bs {u^{N+1}},\bs {\lambda^{N+1}})\in \mathcal{N}_{U}(\bs {u^{N+1}})
\end{equation}
which can be specified by the following equation
\begin{equation}\label{optimal_control_QNLP}
\bs{ u^{N+1}}=\textbf{Pr}_{\lbrack\alpha(t),\beta(t)\rbrack}\left\{\bs {u^{N}}-\frac{g'(\bs {u^{N}})+(\bs{\lambda^{N+1}})^{\top}\bs h_{\bs u}(\bs{x^N},\bs{u^N})}{g''(\bs {u^{N}})}\right\}
\end{equation}

Then (\ref{QLNP_optimal_condtion_1})-(\ref{QLNP_optimality_condition_end}) constitute the first-order optimality conditions for problem $Q^{N+1}$.\\

Remark 4.

1. Note that we can specify the optimal control $\bs{u^{N+1}}$ for every subproblem $Q^{N+1}$ using (\ref{optimal_control_QNLP}). Taking avantage of induction, it is obvious that the sequence $\{(\bs {x^{N}},\bs{ u^{N}})\}$ belongs to $C(\lbrack t_0,t_f\rbrack;\mathbb{R}^n)\times C(\lbrack t_0,t_f\rbrack;\mathbb{R})$ as long as the initial guess $(\bs {x^{0}},\bs{ u^{0}})$ is chosen from functional space $C(\lbrack t_0,t_f\rbrack;\mathbb{R}^n)\times C(\lbrack t_0,t_f\rbrack;\mathbb{R}^n)$.\\

2. The rest convergence analysis in this section is in continuous function space that means the norm is infinity norm, for convinence we directly denote as $\|\cdot\|$ in this section. 

\subsection{Convergence Analysis for Algorithm 1}
\qquad In the following analysis, we prove the convergence for special case  precisely $E=0$ and $\bs e_f=0$ which means the terminal state is unconstrained. Denote the solution sequence obtained by Algorithm 1 as $\{(\bs{x^N},\bs{u^N})\}$, and $\{\bs{\lambda^N}\}$ as the costate sequence each $\bs{\lambda^N}$ is the costate vector for subproblem $Q^{N-1}$ , $\{\bs{\eta^N}\}$ as the multiplier sequence each $\bs{\eta^N}$ is the multiplier corresponding to the terminal state constraint for problem $Q^{N-1}$. Besides, for each $N$ $(\bs{x^N},\bs{u^N},\bs{\lambda^N},\bs{\eta^N})$ is the solution to the optimality conditions of problem $Q^{N-1}$ decribed by (\ref{QLNP_optimal_condtion_1})-(\ref{QLNP_optimality_condition_end}).

We assume that the following hypotheses are satisfied.\\

\textbf{(H1)} Function $f: \mathbb {R}^n\rightarrow\mathbb{R}$ is second-order continuous differentiable, function $g: \mathbb{R}\rightarrow\mathbb{R}$ is third-order continuous differentiable.
The following condition corresponding to $g$ is satisfied 
\begin{equation*}
\sup_{x \in \mathbb{R}}\left|\frac{g'''(x)g'(x)}{(g''(x))^2}\right|<1
\end{equation*}

\textbf{(H2)} Function $\bs h: \mathbb{R}^n\times\mathbb{R}\rightarrow\mathbb{R}^n$ is second-order continuous differentiable corresponding to the arguments and function $\bs \alpha, \bs\beta : \lbrack t_0,t_f\rbrack\rightarrow\mathbb{R}$ is continuous.

\textbf{(H3)} The sequence $\{(\bs{x^N},\bs{u^N})\}$ obtained  by Algorithm 1 is uniformly bounded.

Remark 5.

(1) The condition appears in \textbf{(H1)} can be understood as constraining the curvative of function $g$, if $g$ is a positive-definite quadratic function the condition holds clearly. 

(2) We set $t_0=0$ in the rest proof without loss of generality.

(3) Since the terminal state is unconstrained, we can conclude that $\bs{\eta^N}=\bs0$ for all $N$ which means that at terminal time the costate function $\bs{\lambda^N}$ is zero.

\begin{lemma}
	There exist constant $M_1,M_2,C_1,C_2$ such that for all $N$ the following inequalities hold.
	\begin{align*}
	&|\bs h_{\bs x}(\bs{x^N},\bs{u^N})|\leq M_1\quad
	|\bs h_{\bs u}(\bs{x^N},\bs{u^N})|\leq M_2\quad\\
	&\left|\frac{\partial^2}{\partial x_k\partial x_j} h_i (\bs{x^N},\bs{u^N})\right|\leq C_1\quad \left|\frac{\partial^2}{\partial x_k\partial u} h_i (\bs{x^N},\bs{u^N})\right|\leq C_1 \quad
	\left|\frac{\partial^2}{\partial u^2}h_i \bs(\bs{x^N},\bs{u^N})\right|\leq C_2 \quad \text{for all $i,j,k=1,2,\cdots,n$}
	\end{align*}
\end{lemma}
{\emph{Proof}}. It follows from the hypotheses \textbf{(H2)} and \textbf{(H3)}, the lemma holds obviously.
$\hfill\qedsymbol$\\

For simplicity, we directly use the constant described in the above lemma in the following proof.
We first approximate the distance from $\bs{x^{N+1}}$ to $\bs{x^{N}}$ and the distance from $\bs{\lambda^{N+1}}$ to $\bs{\lambda^{N}}$, which play a key role in analyzing the convergence.

\begin{lemma}\label{distace_x}
	For $\bs{x^{N+1}}$ and $\bs{x^N}$, the following inequality holds for any $t\in\lbrack 0,t_f\rbrack$ 
	\begin{equation*}
	|\bs{x^{N+1}}(t)-\bs{x^{N}}(t)|\leq \int_{0}^t\exp(M_1(t-r))\cdot (2M_1|\bs{x^{N}}(r)-\bs{x^{N-1}}(r)|+M_2|\bs{u^{N+1}}(r)-\bs{u^{N}}(r)|+2M_2|\bs{u^{N}}(r)-\bs{u^{N-1}}(r)|)dr
	\end{equation*}
\end{lemma}
{\emph{Proof}}.
    Because $\bs{x^N}$ is the solution to (\ref{QLNP_dynamic}). Precisely $\bs{x^N}$ and $\bs{x^{N+1}}$ satisfy the following differential equations
	\begin{align*}
	\bs{\dot{{x}}^{N+1}}&=\bs h(\bs {x^{N}},\bs{u^{N}})+
	\bs h_{\bs x}(\bs{ x^{N}},\bs{ u^{N}})(\bs{ x^{N+1}}-\bs{ x^{N}})+
	\bs h_{\bs u}(\bs {x^{N}},\bs {u^{N}})(\bs{u^{N+1}}-\bs {u^{N}})\\
	\bs{\dot{{x}}^{N}}&=\bs h(\bs{ x^{N-1}},\bs{ u^{N-1}})+
	\bs h_{\bs x}(\bs{ x^{N-1}},\bs{ u^{N-1}})(\bs {x^{N}}-\bs {x^{N-1}})+
	\bs h_{\bs u}(\bs{ x^{N-1}},\bs{ u^{N-1}})(\bs {u^{N}}-\bs {u^{N-1}})
	\end{align*}
	
	Subtracting the above two expressions and applying the mean value theorem, we get
	\begin{align*}
	\bs{\dot{{x}}^{N+1}}-\bs{\dot{{x}}^{N}}&=\bs h_{\bs x}(\bs{x^N},\bs{u^N})(\bs{x^{N+1}}-\bs{x^N})+(-\bs h_{\bs x}(\bs{x^{N-1}},\bs{u^{N-1}})+\bs h_{\bs x}(\bs{x^{\xi}},\bs{u^\xi}))(\bs{x^{N}}-\bs{x^{N-1}})\\
	&+\bs h_{\bs u}(\bs{x^N},\bs{u^N})(\bs{u^{N+1}}-\bs{u^N})+(-\bs h_{\bs u}(\bs{x^{N-1}},\bs{u^{N-1}})+\bs h_{\bs u}(\bs{x^{\xi}},\bs{u^\xi}))(\bs{u^{N}}-\bs{u^{N-1}})
	\end{align*}
	Here $(\bs{x^{\xi}},\bs{u^\xi})$ generated by mean value theorem represents $\xi(\bs{x^{N-1}},\bs{u^{N-1}})+(1-\xi)(\bs{x^N},\bs{u^N})$ for some $\xi$.\\
	Since the property of vector modulus and Lemma 7.2, it follows that any time point $t\in\lbrack 0,t_f\rbrack$ there holds the following inequality 
	\begin{equation*}
	|\bs{x^{N+1}}(t)-\bs{x^{N}}(t)|\leq \int_{0}^{t}M_1|\bs{x^{N+1}}(r)-\bs{x^N}(r)|+2M_1|\bs{x^N}(r)-\bs{x^{N-1}}(r)|+M_2|\bs{u^{N+1}}(r)-\bs{u^N}(r)|+2M_2|\bs{u^{N}}(r)-\bs{u^{N-1}}(r)|dr
	\end{equation*}
    By Gronwall's lemma,
	\begin{equation*}
	|\bs{x^{N+1}}(t)-\bs{x^{N}}(t)|\leq \int_{0}^t\exp(M_1(t-r))\cdot (2M_1|\bs{x^{N}}(r)-\bs{x^{N-1}}(r)|+M_2|\bs{u^{N+1}}(r)-\bs{u^{N}}(r)|+2M_2|\bs{u^{N}}(r)-\bs{u^{N-1}}(r)|)dr
	\end{equation*}
	Thus the result is proved. 
	$\hfill\qedsymbol$\\

Next we approximate the distance from $\bs{\lambda^{N+1}}$ to $\bs{\lambda^{N}}$ and the norm of $\bs{\lambda^{N}}$.
\begin{lemma}\label{lambda_norm}
	There exists constant $M_1$, for all costate vector $\bs{\lambda^N}$ there holds that inequality for any $t\in\lbrack 0,t_f\rbrack$
	\begin{align*}
	|\bs{\lambda^N}(t)|\leq \int_t^{t_f}\exp(M_1(t_f-r))\cdot|\nabla f(\bs{x^{N-1}}(r))+\nabla^2 f(\bs{x^{N-1}}(r))(\bs{x^N}(r)-\bs{x^{N-1}}(r))|dr
	\end{align*}	
\end{lemma}
{\emph{Proof}}.
	The proof is similar to the proof of Lemma 7.3, and hence is omitted.  
$\hfill\qedsymbol$
\begin{lemma}\label{distance_lambda}
	There exists contant $M_1,M_2$ and sequence $\{\bs{x^{\xi_N}}\}$ precisely each $\bs{x^{\xi_N}}=\xi_N\bs{x^{N-1}}+(1-\xi_N)\bs{x^{N}}$ (for some $0\leq\xi_N\leq 1$) such that for each $N$, the following inequality holds for any $t\in\lbrack 0,t_f\rbrack$
	\begin{align*}
	|\bs{\lambda^{N+1}}(t)-\bs{\lambda^N}(t)|&\leq \int_{t}^{t_f}(|\nabla^2 f(\bs{x^N}(r))(\bs{x^{N+1}}(r)-\bs{x^N}(r))|
	+|(\nabla^2 f(\bs x^{\bs\xi_N}(r))-\nabla^2f(\bs{x^{N-1}}(r))) (\bs{x^N}(r)-\bs{x^{N-1}}(r))|\\
	&+|(\bs h_{\bs x}(\bs{x^N}(r),\bs{u^N}(r))
	-\bs h_{\bs{x}}(\bs{x^{N-1}}(r),\bs{u^{N-1}}(r)))\bs{\lambda^N}(r)|)\exp(M_1({t_f-r}))
	dr
	\end{align*}
\end{lemma}
{\emph{Proof}}.
	Since $\bs{\lambda^N}$ respresent costate vector, it follows from (\ref{QLNP_optimal_condtion_1}) and (\ref{QLNP_con_cos})
	\begin{align*}
	\bs{\dot{{\lambda}}^{N+1}}&=-\bs h_{\bs x}(\bs{ x^{N}},\bs{ u^{N}})^{\top}\bs{\lambda^{N+1}}
	-\nabla f(\bs{x^{N}})-\nabla^2 f(\bs{x^{N}})(\bs{x^{N+1}-x^{N}})\\
	\bs{\dot{{\lambda}}^{N}}&=-\bs h_{\bs x}(\bs{ x^{N-1}},\bs{ u^{N-1}})^{\top}\bs{\lambda^{N}}
	-\nabla f(\bs{x^{N-1}})-\nabla^2 f(\bs{x^{N-1}})(\bs{x^{N}-x^{N-1}})\\
	\bs{\lambda^{N+1}}(t_f)&=0\\
	\bs{\lambda^{N}}(t_f)&=0
	\end{align*}
	
    The rest proof is similar to the proof of Lemma 7.3, and hence is omitted.
	$\hfill\qedsymbol$\\

With these preparations, we are now able to prove the convergence of alogorithm 1 for solving nonlinear problem (NP) without terminal state constraint.
\begin{theorem}\label{conver1}
There exists some constant $T>0$, as long as $t_f<T$ there exists $(\bs {x^{\infty}},\bs {u^{\infty}})$ in functional space $C(\lbrack t_0,t_f\rbrack;\mathbb{R}^n)\times C\lbrack t_0,t_f\rbrack$ such that sequence $\{(\bs {x^{N}},\bs {u^{N}})\}_{N=1}^{\infty}$ converge to $(\bs {x^{\infty}},\bs {u^{\infty}})$
\end{theorem}
{\emph{Proof}}.
Firstly we can specify $\bs{u^N}$ by (\ref{optimal_control_QNLP}), by the property of projection operator one can derive that
\begin{align*}
&\|\bs{u^{N+1}}-\bs{u^N}\|=\left\|\bs{u^N}-\frac{g'(\bs{u^N})+(\bs{\lambda^{N+1}})^{\top}\bs{h_{u}(x^N,u^N)}}{g''(\bs{u^N})}-\bs{u^{N-1}}+\frac{g'(\bs{u^{N-1}})+(\bs{\lambda^{N}})^{\top}\bs{h_{u}(x^{N-1},u^{N-1})}}{g''(\bs{u^{N-1}})}\right\|\\
&\leq \left\|(\bs{u^{N}}-\bs{u^{N-1}})-\left(\frac{g'(\bs{u^N})}{g''(\bs{u^N})}-\frac{g'(\bs{u^{N-1}})}{g''(\bs{u^{N-1}})}
\right)\right\|
+\left\|\frac{(\bs{\lambda^{N+1}})^{\top}\bs{h_{u}(x^N,u^N)}}{g''(\bs{u^N})}-\frac{(\bs{\lambda^{N}})^{\top}\bs{h_{u}(x^{N-1},u^{N-1})}}{g''(\bs{u^{N-1}})}\right\|	\\
&\leq\left\|\frac{g'(\bs{u^{\xi}})g'''(\bs{u^{\xi}})}{(g''(\bs{u^{\xi}}))^2}(\bs{u^N}-\bs{u^{N-1}})\right\|+\left\|\frac{(\bs{\lambda^{N+1}-\lambda^{N}})^{\top}\bs{h_{u}(x^N,u^N)}}{g''(\bs{u^N})}\right\|\\
&+\left\|(\bs{\lambda^N})^{\top}(\frac{\bs{h_{u}(x^N,u^N)}}{g''(\bs{u^N})}-\frac{\bs{h_{u}(x^{N-1},u^{N-1})}}{g''(\bs{u^{N-1}})})\right\|
\end{align*}

Taking advantage of hypothese \textbf{(H1)}, we know that there exist some constant $0<\theta<1$ such that $\left|\frac{g'''(\bs{u^{\xi}})g'(\bs{u^{\xi}})}{(g''(\bs{u^{\xi}}))^2}\right|<\theta$.

From Lemma 7.2 and assumption \textbf{(H1)} and \textbf{(H3)}, by mean value theorem, one can derive that there exist some constant $\widehat{M},\widetilde{M_1},\widetilde{M_2}$ such that
\begin{align*}
\left\|\frac{\bs{h_{u}(x^N,u^N)}}{g''(\bs{u^N})}\right\|\leq\widehat{M},  \left\|\frac{\bs{h_{u}(x^N,u^N)}}{g''(\bs{u^N})}-\frac{\bs{h_{u}(x^{N-1},u^{N-1})}}{g''(\bs{u^{N-1}})}\right\|\leq \widetilde{M_1}\|\bs{x^{N}}-\bs{x^{N-1}}\|+\widetilde{M_2}\|\bs{u^{N}}-\bs{u^{N-1}}\|
\end{align*}

Thus the following inequality holds
\begin{equation*}
\begin{aligned}
\|\bs{u^{N+1}}-\bs{u^N}\|&\leq\theta\|\bs{u^N}-\bs{u^{N-1}}\|+ \widehat{M}\|\bs{\lambda^{N+1}-\lambda^{N}}\|+\|\bs{\lambda^N}\|(\widetilde{M_1}\|\bs{x^{N}}-\bs{x^{N-1}}\|+\widetilde{M_2}\|\bs{u^{N}}-\bs{u^{N-1}}\|) 
\end{aligned}
\end{equation*}

Since we have approximate $|\bs{\lambda^N}(t)|$ by Lemma 7.4 and $|\bs{\lambda^{N+1}}(t)-\bs{\lambda^N}(t)|$ by Lemma 7.5, combining with the assumption \textbf{(H1)}, \textbf{(H2)}, it follows that there exists constant $C_1,C_2,C_3,\widetilde{C_3},\widetilde{C_4}$ such that

\begin{equation}\label{eq_1}
\begin{aligned}
\|\bs{u^{N+1}}-\bs{u^{N}}\|&\leq \theta\|\bs{u^N}-\bs{u^{N-1}}\|+(\exp(M_1t_f)-1)(\widetilde{C_3}\|\bs{x^{N}}-\bs{x^{N-1}}\|+\widetilde{C_4}\|\bs{u^{N}}-\bs{u^{N-1}}\|)\\
	&+\frac{\exp(M_1t_f)-1} {M_1}(C_1\|\bs{x^{N+1}}-\bs{x^N}\|+C_2\|\bs{x^{N}}-\bs{x^{N-1}}\|+C_3\|\bs{u^{N}}-\bs{u^{N-1}}\|)
\end{aligned}
\end{equation}

Substituting $t$ in right side of the inequality described in Lemma 7.3 by $t_f$, one can derive that
\begin{equation}\label{eq_2}
\|\bs{x^{N+1}}-\bs{x^N}\|\leq\frac{\exp(M_1t_f)-1}{M_1}(2M_1\|\bs{x^{N}}-\bs{x^{N-1}}\|+M_2\|\bs{u^{N+1}}-\bs{u^{N}}\|+2M_2\|\bs{u^{N}}-\bs{u^{N-1}}\|)
\end{equation}

Combining inequalities (\ref{eq_1}) and (\ref{eq_2}), the following inequality holds 
\begin{equation}\label{ineq_3}
\begin{aligned}
&(1-\frac{\exp(M_1t_f)-1}{M_1}C_1)\|\bs{x^{N+1}}-\bs{x^N}\|+(1-\frac{\exp(M_1t_f)-1}{M_1}M_2)\|\bs{u^{N+1}}-\bs{u^{N}}\|\\
\leq&(\frac{\exp(M_1t_f)-1} {M_1}C_2+(\exp(M_1t_f)-1)\widetilde{C_3}+2(\exp(M_1t_f)-1))\|\bs{x^{N}}-\bs{x^{N-1}}\|\\
&+(\theta+\frac{\exp(M_1t_f)-1} {M_1}C_3+(\exp(M_1t_f)-1)\widetilde{C_4}+2\frac{M_2}{M_1}(\exp(M_1t_f)-1))\|\bs{u^{N}}-\bs{u^{N-1}}\|
\end{aligned}
\end{equation}

Considering the following inequalities
\begin{equation}\label{ineq_1}
\left\{	
\begin{aligned}
1&>&\theta +\frac{\exp(M_1t_f)-1} {M_1}C_3+(\exp(M_1t_f)-1)\widetilde{C_4} +\frac{3M_2}{M_1}(\exp(M_1t_f)-1)\notag\\
1&>&\frac{\exp(M_1t_f)-1} {M_1}(C_1+C_2)+2(\exp(M_1t_f)-1)+(\exp(M_1t_f)-1)\widetilde{C_3}\notag
\end{aligned}
\right.	
\end{equation}
Beacause the right hand side is continuous corresponding to $t_f$, we conclude that there exists a constant $T$ such that the above inequalities holds for all $t_f\in\lbrack 0,T\rbrack$.

If we choose $t_f(t_f\in\lbrack 0,T\rbrack)$, by (\ref{ineq_1}) we can directly obtain that there exists some positive constant $\gamma<1$ such that
\begin{equation}
\begin{aligned}
\|\bs{x^{N+1}}-\bs{x^N}\|+\|\bs{u^{N+1}}-\bs{u^N}\|&\leq
\gamma(\|\bs{x^{N}}-\bs{x^{N-1}}\|+\|\bs{u^{N}}-\bs{u^{N-1}}\|)
\end{aligned}
\end{equation}
	
	Then one have that for arbitrarily $N1>N2$ there holds
	\begin{equation*}
	\begin{aligned}
	\|\bs{x^{N1}}-\bs{x^{N2}}\|+\|\bs{u^{N1}}-\bs{u^{N2}}\|&\leq
	\sum_{K=N2}^{N1-1} (\|\bs{x^{K+1}}-\bs{x^{K}}\|+\|\bs{u^{K+1}}-\bs{u^{K}}\|)\\
	&\leq\sum_{K=N2}^{N1-1}\gamma^{K-1}(\|\bs{x^{2}}-\bs{x^{1}}\|+\|\bs{u^{2}}-\bs{u^{1}}\|)
	\end{aligned}
	\end{equation*}
	
	Thus we know that the sequence $\{(\bs {x^{N}},\bs {u^{N}})\}_{N=1}^{\infty}$ is a Cauthy sequence in functional space $C(\lbrack t_0,t_f\rbrack;\mathbb{R}^n)\times C\lbrack t_0,t_f\rbrack$. 
	Because the functional space considered is a Banach space thus the therorem holds explicitly.
	
	Hence the theorem is proved. 
	$\hfill\qedsymbol$
    \begin{theorem}
    	The limit point $(\bs{x^\infty},\bs{u^\infty})$ obtained by the sequence $\{(\bs{x^N},\bs{u^N})\}_{N=1}^{\infty}$ described in Theorem 7.1 satisfies the necessary conditions for original problem (NP) (\ref{QLNP_optimal_condtion_1})-(\ref{QLNP_optimality_condition_end}).
    \end{theorem}
    {\emph{Proof}}.
     	Review the optimality conditions for subproblem $Q^{N+1}$ (\ref{QLNP_optimal_condtion_1})-(\ref{QLNP_optimality_condition_end}), for sequence $\{(\bs{x^N},\bs{u^N})\}_{N=1}^{\infty}$ we have	
     	\begin{align*}
     	&\bs{{\dot{ x}}^{N+1}}=\bs h(\bs {x^{N}},\bs {u^{N}})+
     	\bs h_{\bs x}(\bs{ x^{N}},\bs {u^{N}})(\bs {x^{N+1}}-\bs {x^{N}})+
     	\bs h_{\bs u}(\bs{ x^{N}},\bs {u^{N}})(\bs {u^{N+1}}-\bs {u^{N}})\\
     	&\bs{{\dot{\lambda}}^{N+1}}=-\nabla^2f(\bs{ x^{N}})(\bs {x^{N+1}}-\bs{ x^{N}})-\nabla f(\bs {x^{N}})-\bs h_{\bs x} (\bs {x^{N}},\bs{ u^{N}})^{\top}\bs{\lambda^{N+1}} \\
     	&\bs {x^{N+1}}(t_0)=\bs x_{0}\\
     	&\bs {\lambda^{N+1}}(t_f)=0\\
     	&\bs -H _{\bs u}(\bs{ x^{N+1}},\bs{ u^{N+1}},\bs {\lambda^{N+1}})\in \mathcal{N}_{U}(\bs {u^{N+1}})
     	\end{align*}
     	As is proved in Theorem 7.1, sequence $\{(\bs {x^{N}},\bs {u^{N}})\}_{N=1}^{\infty}$ converges to $(\bs {x^{\infty}},\bs {u^{\infty}})$ and 	$\lim\limits_{k\to\infty}\|\bs{x^{N+1}}-\bs{x^{N}}\|+\|\bs{u^{N+1}}-\bs{u^{N}}\|=0$.
     	
     	Thus we obtain the following equations by taking limits upon $N$ to infinity both sides
     	\begin{align*}
     	&\bs{{\dot{ x}}^{\infty}}=\bs h(\bs {x^{\infty}},\bs {u^{\infty}})\\
     	&\bs {x^{\infty}}(t_0)=\bs x_{0}
     	\end{align*}
     	
       Since the costate vector is solution to linear differential equations, we can represent $\bs{\lambda^{N+1}}$ by the following equations
       \begin{equation*}
       \bs {\lambda^{N+1}}=\int_{t}^{t_f} \Phi^N(t,r)(\nabla^2f(\bs{ x^{N}})(\bs {x^{N+1}}-\bs{ x^{N}})+\nabla f(\bs {x^{N}}))dr
       \end{equation*}
       $\Phi^N(t,r)=\Phi^N(t)\Phi^N(r)^{-1}$. $\Phi^N(t)$ denotes the solution to the following differential equations where $I_n$ denotes an n by n identity matrix.  
       \begin{equation*}
       \dot{\Phi}(t)=\bs h_{\bs x}(\bs {x^{N}},\bs{ u^{N}})^{\top}\Phi(t),\qquad \Phi(0)=I_{n}
       \end{equation*}
       
       It is shown in \cite{MR319376} that $\Phi^N(t)$ tends to $\Phi^\infty(t)$ when $N$ tends to infinity. Here $\Phi^\infty(t)$ is solution to the following differential equations.
       \begin{equation*}
       \dot{\Phi}(t)=\bs h_{\bs x}(\bs {x^{\infty}},\bs{ u^{\infty}})^{\top}\Phi(t),\qquad \Phi(0)=I_{n}
       \end{equation*}
       
       Thus it holds that $\bs {\lambda^{N}}$ tends to $\bs {\lambda^{\infty}}$ while $N$ tends to infinity.
       \begin{equation*}
       \bs {\lambda^{\infty}}=\int_{t}^{t_f} \Phi^\infty(t,r)\nabla f(\bs {x^{\infty}})dr
       \end{equation*}
      
      We conclude that $\bs {\lambda^{\infty}}$ satisfy
      \begin{equation*}
      \begin{aligned}
      &\bs{{\dot{ \lambda}}^{\infty}}=-\bs h_{\bs x}(\bs {x^{\infty}},\bs {u^{\infty}})^{\top}\bs{\lambda^\infty}-\nabla f(\bs{x^\infty})\\
      &\bs {\lambda^{\infty}}(t_f)=\bs 0
      \end{aligned}
      \end{equation*}
          
       Because of the definition of normal cone, we can directly obtain 
    \begin{equation*}
   \bs -H_{\bs u}(\bs {x^{\infty}},\bs {u^{\infty}},\bs {\lambda^{\infty}})\in \mathcal{N}_{U}(\bs {u^{\infty}})
    \end{equation*}
   
Thus we complete the proof
$\hfill\qedsymbol$

Remark 6.

We note that for terminal state unconstrained case the Algorithm 1 is of globally convergence in some time horizon and the convergence rate is linear which is fast. But we should also notice that these conclusion is holds while time horizon satisfy some condition.  

\subsection{Convergence Analysis for Algorithm 2}

We will prove the convergence property of Algorithm 2 for general case that the terminal state constraints cannot be omit.

We denote sequence $\{(\bs{\bar{x}^N},\bs{\bar{u}^N})\}$, $\{\bs{\bar{\lambda}^N}\}$, $\{\bs{\eta^N}\}$, each $(\bs{\bar{x}^N},\bs{\bar{u}^N})$ is obained by solving problem $Q^{N-1}$  and $\{\bs{\bar{\lambda}^N}\}$ is the corresponding costate vector $\bs{\eta^N}$ is the multiplier with regard to terminal state constraints. Besides, for each $N$, $(\bs{\bar{x}^N},\bs{\bar{u}^N},\bs{\bar{\lambda}^N},\bs{\eta^N})$ is the solution to the optimality conditions of problem $Q^{N-1}$ described by (\ref{QLNP_optimal_condtion_1})-(\ref{QLNP_optimality_condition_end}).
 
We assume that following conditions holds in our proof.
 
\textbf{(A1)} Function $f: \mathbb {R}^n\rightarrow\mathbb{R}$ is second-order continuous differentiable, function $g: \mathbb{R}\rightarrow\mathbb{R}$ is second-order continuous differentiable.

\textbf{(A2)} Function $\bs h: \mathbb{R}^n\times\mathbb{R}\rightarrow\mathbb{R}^n$ is second-order continuous differentiable corresponding to the arguments and function $\bs \alpha, \bs\beta : \lbrack t_0,t_f\rbrack\rightarrow\mathbb{R}$ is continuous.

\textbf{(A3)} The obtained sequence $\{(\bs{\bar{x}^N},\bs{\bar{u}^N})\}$ is uniformly bounded and the multiplier sequence $\{\bs{\eta^N}\}$ is bounded.

\textbf{(A4)} The control sequence $\{\bs{\bar{u}^N}\}$ is uniformly equicontinuous.
%\textbf{(A4)} The left and right derivative of $\bs{\bar{u}^N}$ exist for all $N$, and are uniformly bounded.

%Remark 7.

%Assumption \textbf{(A4)} means that the control sequence is nonsingular precisely bang-bang control situtation doesn't appear in our analysis.\\

Review the merit functional $P(\bs x,\bs u)$ introduced in our Algorithm 2.
\begin{equation*}
P(\bs x,\bs u)=\int_{t_0}^{t_f} f(\bs x)+g(\bs u) dt+\theta\|\dot{\bs x}-\bs h(\bs x,\bs u)\|_1
\end{equation*}

We should firstly prove that the line-search is reasonable, thus alogorithm 2 is executable.

\begin{theorem}
	Since $(\bs{\bar{x}^{N+1}},\bs{\bar{u}^{N+1}})$ is the solution to problem $Q^{N+1}$ and the corresponding costate vector $\bs{\bar{\lambda}^{N+1}}$, if there holds $\|\bs{\bar{\lambda}^{N+1}}\|<\theta$ then the derivative of function $P_1(\kappa):=P((\bs {x^N},\bs {u^N})+\kappa\bs {d^N})$ at $\kappa=0$ is negative which guarantee the resonability of line search. Furthermore we can conclude that the sequence $\{P(\bs{ x^N},\bs{u^N})\}$ is monotonically decreasing.
	
	Here $\bs {d^N}:=(\bs{\bar{x}^{N+1}},\bs{\bar{u}^{N+1}})-(\bs{x^{N}},\bs{u^{N}})$, $\bs {d_x^N}=\bs{\bar{x}^{N+1}}-\bs{x^{N}}$, $\bs {d_u^N}=\bs{\bar{u}^{N+1}}-\bs{u^{N}}$.
\end{theorem}
{\emph{Proof}}. We refer to the LEMMA 4.2 . Since it follows similarly without much more effort, we omit
proof here.
$\hfill\qedsymbol$

Then we prove that sequence $\{(\bs{x^N},\bs{u^N})\}$ exists cluster points and any cluster point satisfy the optimality conditions of (NP).

\begin{lemma}
	The sequence $\{(\bs{x^N},\bs{u^N})\}$ obtained is campact in functional space $C(\lbrack t_0,t_f\rbrack;\mathbb{R}^n)\times C(\lbrack t_0,t_f\rbrack;\mathbb{R})$ which means that any subsequence of $\{(\bs{x^N},\bs{u^N})\}$ has cluster point in functional space $C(\lbrack t_0,t_f\rbrack;\mathbb{R}^n)\times C(\lbrack t_0,t_f\rbrack;\mathbb{R})$
\end{lemma}

\emph{Proof}. From assumption \textbf{(A3)} and $(\bs{x^N},\bs{x^N})$ obtained by Algorithm 2, we know that the sequence $\{(\bs{x^N},\bs{u^N})\}$ is uniformly bounded by induction. Taking advantage of assumption \textbf{(A3)}, there holds that the derivative of $\bs{x^N}$ is uniformly bounded. By the equicontinuity of functional sequence discussed in \cite{MR1476423}, we conclude that $\{\bs{x^N}\}$ is equicontinuous. According to the Arzela-Ascoli Theorem shown in \cite{MR1476423}, the sequence $\{(\bs{x^N},\bs{u^N})\}$ is campact in functional space $C(\lbrack t_0,t_f\rbrack;\mathbb{R}^n)\times C(\lbrack t_0,t_f\rbrack;\mathbb{R})$. Thus we complete the proof.
$\hfill\qedsymbol$

\begin{theorem}
	The cluster point of sequence $\{(\bs{x^N},\bs{u^N})\}$ satisfy the optimality conditions of problem (NP)
\end{theorem}
{\emph{Proof}}. We refer to the Theorem 4.3 . Since it follows similarly without much more effort, we omit
proof here.
$\hfill\qedsymbol$
\section{Numerical Experiments}
\qquad In this section, we carry out numerical experiments using five example, where in each problem the control has simple lower and upper bounds.

In each example we solve the subproblem $Q^{N+1}$ obtained by quasilinearization technique by Euler discretization scheme. Dontchev, Hager and Malanowski in \cite{MR1779739} present a convergence result for Euler discretization of control-constrained optimal control problem. 

We called our prensent method as sequential dual method since we solve the dual problem of subproblem $Q^{N+1}$ camparing to solve the subproblem $Q^{N+1}$ directly which is called as sequential primal method. 

We use IPOPT version 3.12.3 \cite{MR2899152} to solve the dicretized optimization problem which is described by a large-scale quadratic programming. We peform all computations on a computer with the process, Inter(R) Core(TM) i5-5200U CPU at 2.20GHz, and with a 8.00-GB RAM.

We set stopping tolerence $tol$=1e-5 in our present algorithm. Besides Ipopt parameter is set by max\_iter=2000, tol=1e-8 and acceptable\_tol=1e-12 for solving each subproblem $Q^{N+1}$ by Euler discretization.

Example 1. Consider the following control-unconstrained optimal control problem

\begin{align*}
\min \quad &\frac{1}{2}\int_{0}^{5} x_1^{2}+x_2^2+u^2 dt\\
\text{subject to}\quad &\dot{x_1}=-x_1+x_2\\
&\dot{x_2}=-0.5x_1-0.5x_2(1-(2+cos(2x_1))^2)+(2+cos(2x_1))u\\
& x(0)=(\frac{\pi}{3}\quad \frac{\pi}{4})^{\top} 
\end{align*}

The nonlinear problem above is of the form (NP) with $E=0$ and $\bs e_f=0$ that terminal state is unconstrained.

Table 1 lists the numerical results for this example with various $N$ of discretization subintervals. As expected, by our discussion on the convergence of Algorithm 1 in section 7, strong duality of linear-quadratic subproblem in section 5 and convergence of Euler discretization given in \cite{MR1779739}, as the partition becomes finer, the optimal value given by the present algorithm called sequential dual method tends to some constant that is identical to the other two methods.

\begin{table}[H]
\caption{Example 1 numerical performance }
\begin{tabular}{| l ||c|c | c| c | c| c| c | c |}
	\hline 
	\multirow{2}{*}{N} & \multicolumn{3}{|c|}{Optimal Value} & \multicolumn{2}{|c|}{number of quasilinearization} & \multicolumn{3}{|c|}{IPOPT time[sec]} \\
	\cline{2-9}
	\multirow{2}{*}{~} &Euler & SeqPrimal & SeqDual   & SeqPrimal & SeqDual & Euler & SeqPrimal & SeqDual \\
	\hline
	\hline
	50  &0.6426  &0.6426 &0.5234 &6 &3 &0.010 &0.034  &0.011 \\
	\hline
	100  &0.6121  &0.6121 &0.5532 &6 &4 &0.014 &0.037  &0.014 \\
	\hline
	200  &0.5972  &0.5972 &0.5680 &6 &4 &0.020 &0.052  &0.020 \\
	\hline
	500  &0.5884 &0.5884 &0.5767    & 6  &4   &0.045 & 0.083 & 0.032 \\
	\hline
	1000 &0.5855  &0.5855  &0.5797     & 6 &4  & 0.108  & 0.119  & 0.051  \\
	\hline
	2000 &0.5840  &0.5840  &0.5811     & 6 &4  & 0.128  & 0.339  & 0.089  \\
	\hline
	5000 &0.5832 & 0.5832  & 0.5820    & 6  & 4 &0.382  & 0.644  &0.308    \\
	\hline
	10000 &0.5829   & 0.5829 & 0.5829   & 6  &4   & 0.716 &2.236 & 0.534  \\
	\hline
\end{tabular}
\end{table}

\begin{figure}[htbp]
	\centering
	\begin{minipage}[t]{0.5\textwidth}
		\centering
		\includegraphics[width=6cm]{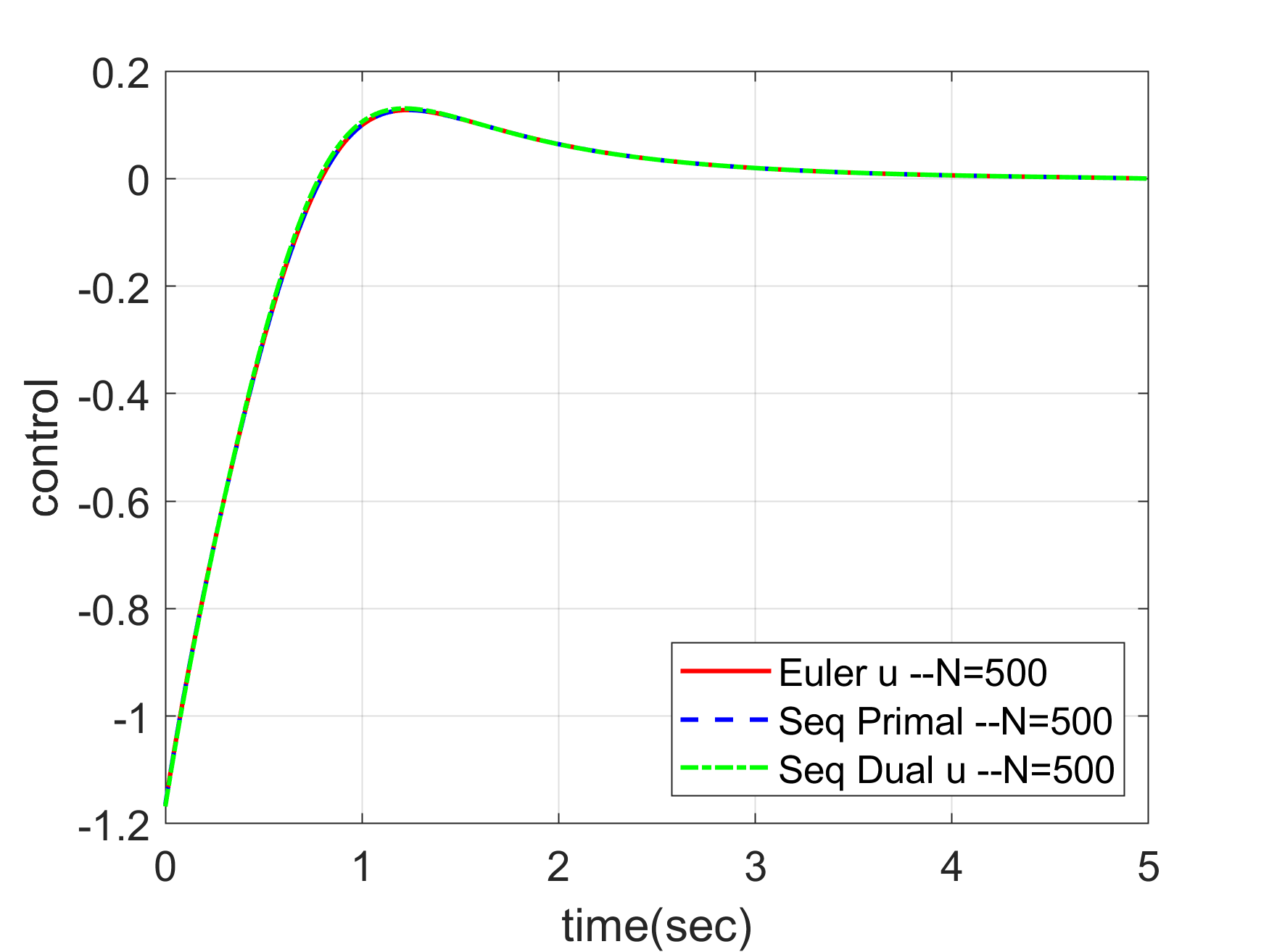}
		%	\caption{N=500时对偶方法的解与直接欧拉离散的解}
	\end{minipage}
	\begin{minipage}[t]{0.49\textwidth}
		\centering
		\includegraphics[width=6cm]{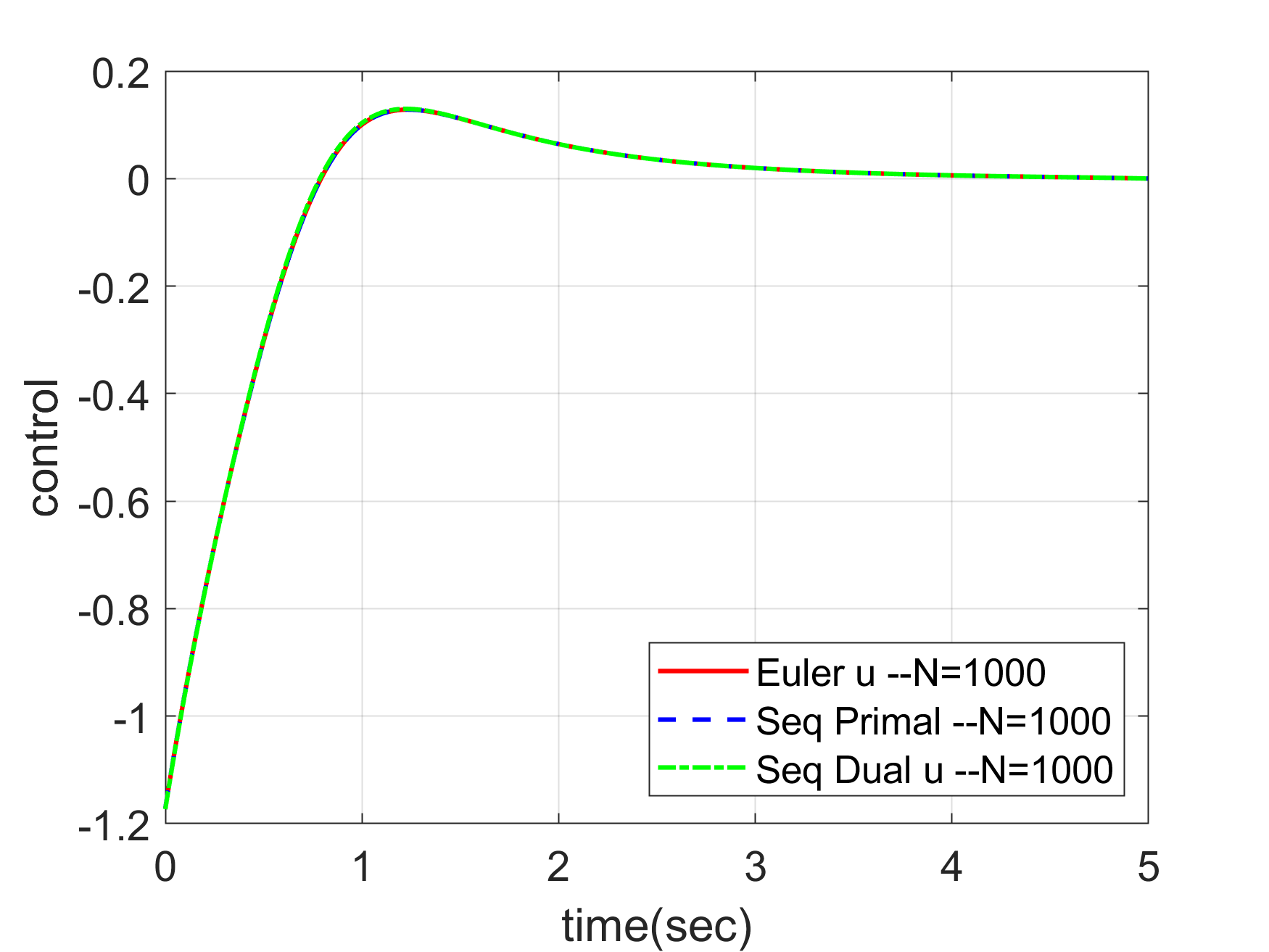}
		%	\caption{N=1000时对偶方法的解与直接欧拉离散的解}
	\end{minipage}
	\caption{Example 1 optimal control obtained with N=500 and N=1000}
\end{figure}

Table 1 also tabulates the quasilinearzation number, note that total quasilinearization number is almost fixed as the partition becoming finer which verify our discussion in section 7. The Ipopt time for solving Example 1 by different method is interesting, note that the objective functional of Example 1 is strongly convex as discussed in Remark 3 indeed the dual problem of subproblem has only the terminal state constrained which can be used to interpret the reason for the IPOPT time camparsion.

Figure 1 depicts the optimal control obtained by solving problem (NP) with $N=500$ and $N=1000$ that nicely illustrates that $N$ get larger squential dual method, sequetial primal method and Euler method approach the continuous-time solution consistently. The convergence result obtained in section 7 is verified. 

Example 2. The following optimal control problem comes from a continuous stirred tank reactor system was studied in \cite{dynamic_1}.
\begin{align*}
\min \quad &\int_{0}^{0.78} x_1^{2}+x_2^2+0.1\cdot u^2 dt\\
\text{subject to}\quad &\dot{x_1}=-2(x_1+0.25)+(x_2+0.5)\exp(\frac{25x_1}{x_1+2})-(x_1+0.25)\cdot u\\
&\dot{x_2}=0.5-x_2-(x_2+0.5)\exp(\frac{25x_1}{x_1+2})\\
&-1\leq u\leq 1\\
& x(0)=(0.05\quad 0)^{\top} 
\end{align*}
The nonlinear problem above is of the form (NP) with $E=0$ and $\bs e_f=0$ that the terminal state is unconstrained. Besides, the objective functional of this problem is also strongly convex in state variables and control variable.

Table 2 lists the numerical results for this example that are generated in the same way as Table 1. The results can be interpreted similarly. Besides, our present method also behaves well when the partition is finer.
\begin{table}[H]
	\caption{Example 2 numerical performance }
\begin{tabular}{|  l ||c|c | c| c| c| c| c | c  |}
	\hline 
	\multirow{2}{*}{N} & \multicolumn{3}{|c|}{Optimal Value} & \multicolumn{2}{|c|}{number of quasilinearization} & \multicolumn{3}{|c|}{IPOPT time[sec]} \\
	\cline{2-9}
	\multirow{2}{*}{~} &Eluer & SeqPrimal & SeqDual & SeqPrimal & SeqDual & Euler & SeqPrimal & SeqDual \\
	\cline{2-9}
	\hline
	\hline
	50 & 0.0301 &0.0301 &0.0139 &11 &10 &0.047&0.278 &0.039 \\
	\hline
	100 & 0.0295 &0.0295 &0.0286 &11 &11 &0.058&0.336 &0.076 \\
	\hline
	200 & 0.0292 &0.0292 &0.0286 &11 &11 &0.131&0.482 &0.076 \\
	\hline
	500 & 0.0291&0.0291  & 0.0288  &11& 11&0.745&0.729 &0.126  \\
	\hline
	1000 & 0.0290 &0.0290 &0.0288  &11 &11 &1.21&1.09 &0.218 \\
	\hline
	2000 & 0.0290 &0.0292 &0.0286 &11 &11 &3.62&3.17 &0.342 \\
	\hline
	5000 & 0.0290 &0.0290 &0.0290 &11 &11 &4.31& 5.16 &0.948 \\
	\hline
	10000 & 0.0290 &0.0290 &0.0290  &11 &11 &9.68 &9.84 &2.436 \\
	\hline
\end{tabular}
\end{table}
\begin{figure}[htbp]
	\centering
	\begin{minipage}[t]{0.49\textwidth}
		\centering
		\includegraphics[width=6cm]{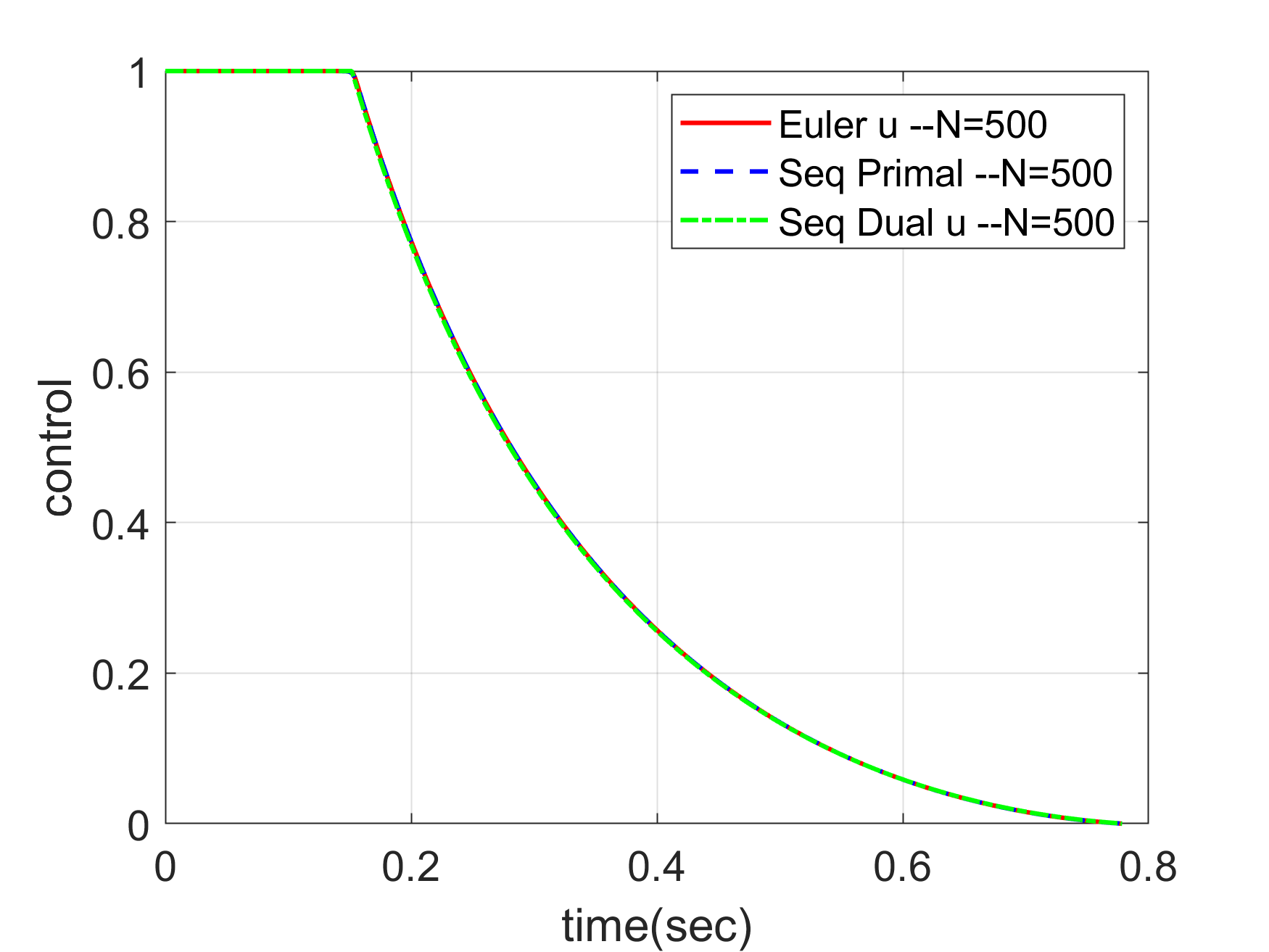}
		%	\caption{N=500时对偶方法的解与直接欧拉离散的解}
	\end{minipage}
	\begin{minipage}[t]{0.49\textwidth}
		\centering
		\includegraphics[width=6cm]{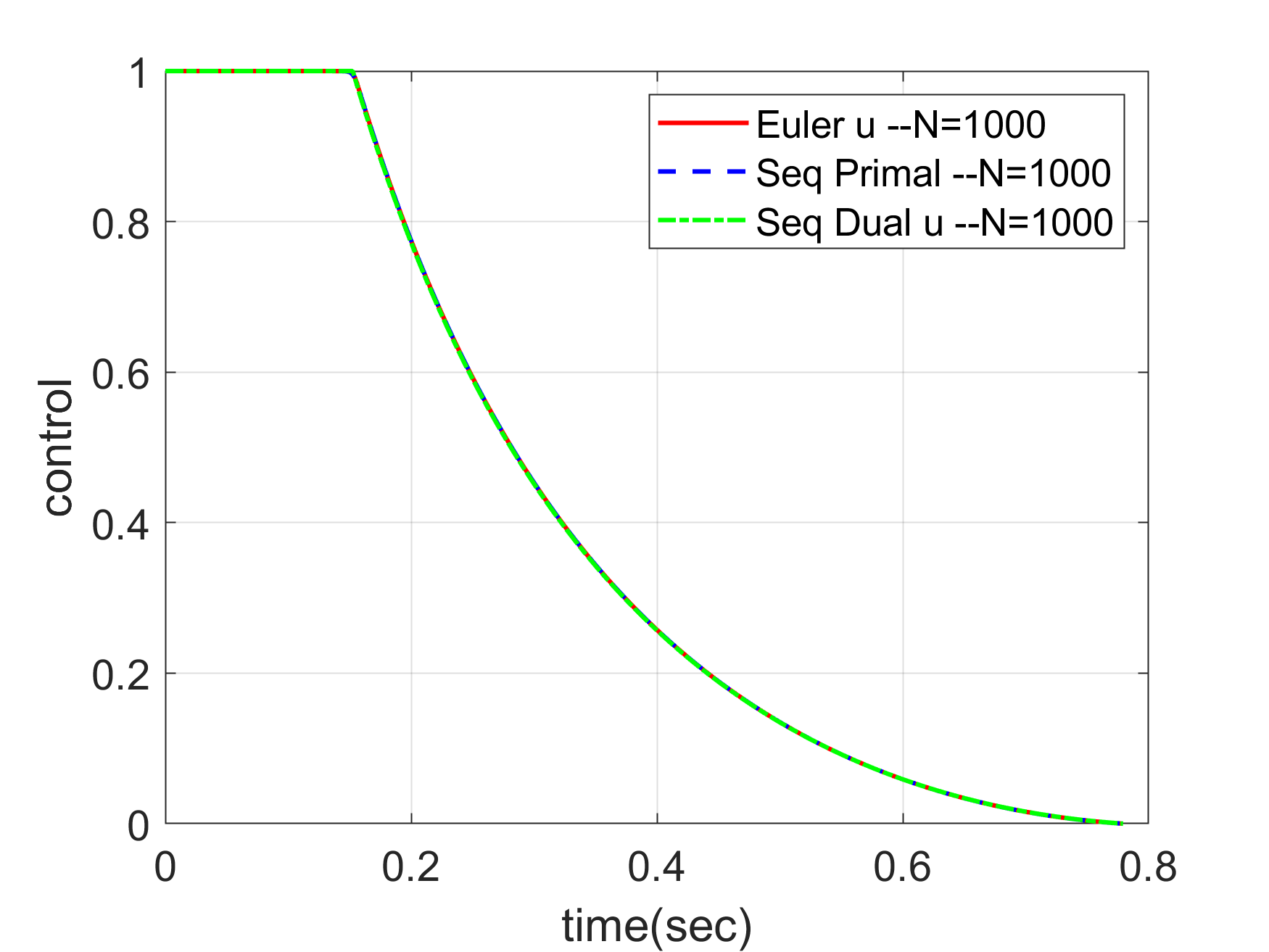}
		%	\caption{N=1000时对偶方法的解与直接欧拉离散的解}
	\end{minipage}
	\caption{Example 2 optimal control obtained with N=500 and N=1000}
\end{figure}

Example 3.  The following optimal control problem which is called Rayleigh problem \cite{MR1860620} is taken into consideration.
\begin{align*}
\min \quad &\frac{1}{2}\int_{0}^{4.5} x_1^{2}+u^2 dt\\
\text{subject to}\quad &\dot{x_1}=x_2\\
&\dot{x_2}=(1.4-0.14x_2^2)\cdot x_2-x_1+4u\\
&-1\leq u\leq 1 \\
& x(0)=(-5\quad -5)^{\top} 
\end{align*}
The nonlinear problem above is of the form (NP) with $E=0$ and $\bs e_f=0$, the following table shows the computation results. Besides, the objective functional of this problem is only convex corresponding to the state variables.

Table 3 lists the numerical results for this example that are generated in the same way as Table 1. 
Notice that our present method also behaves well when the partition is finer.
\begin{table}[H]
	\caption{Example 3 numerical performance }
\begin{tabular}{| l || c | c | c| c| c | c|c |c|}
	\hline 
	\multirow{2}{*}{N} & \multicolumn{3}{|c|}{Optimal Value} & \multicolumn{2}{|c|}{number of quasilinearization} & \multicolumn{3}{|c|}{IPOPT time[sec]} \\
	\cline{2-9}
	\multirow{2}{*}{~} &Eluer & SeqPrimal & SeqDual  & SeqPrimal & SeqDual & Euler & SeqPrimal & SeqDual \\
	
	\hline
	\hline
	50  & 23.2885&23.2425 & 21.2096 &12 & 12 &0.031&0.282 &0.079\\
	\hline
	100  & 22.9310&22.9082 & 21.6877 &12 & 11 &0.073&0.347 &0.107\\
	\hline
	200  & 22.6410&22.6297 & 21.9850 &12 & 11 &0.121&0.497 &0.147\\
	\hline
	500  & 22.4458&22.4413 & 22.1763  &12 &11  &0.448&0.813 &0.278 \\
	\hline
	1000 & 22.3780&22.3757 & 22.2421  &12 &11&0.779&1.43 &0.433 \\
	\hline
	2000  & 22.3436&22.3425 & 22.2754 &12 & 11 &1.13&2.295 &0.7144\\
	\hline
	5000 & 22.3228&22.3228 & 22.2955&12 &11 &2.81&5.47 &1.913  \\
	\hline
	10000 & 22.3159&22.3159 & 22.3121 &12 & 11 &6.1& 13.66 &4.973 \\
	\hline
\end{tabular}
\end{table}
\begin{figure}[htbp]
	\centering
	\begin{minipage}[t]{0.49\textwidth}
		\centering
		\includegraphics[width=6cm]{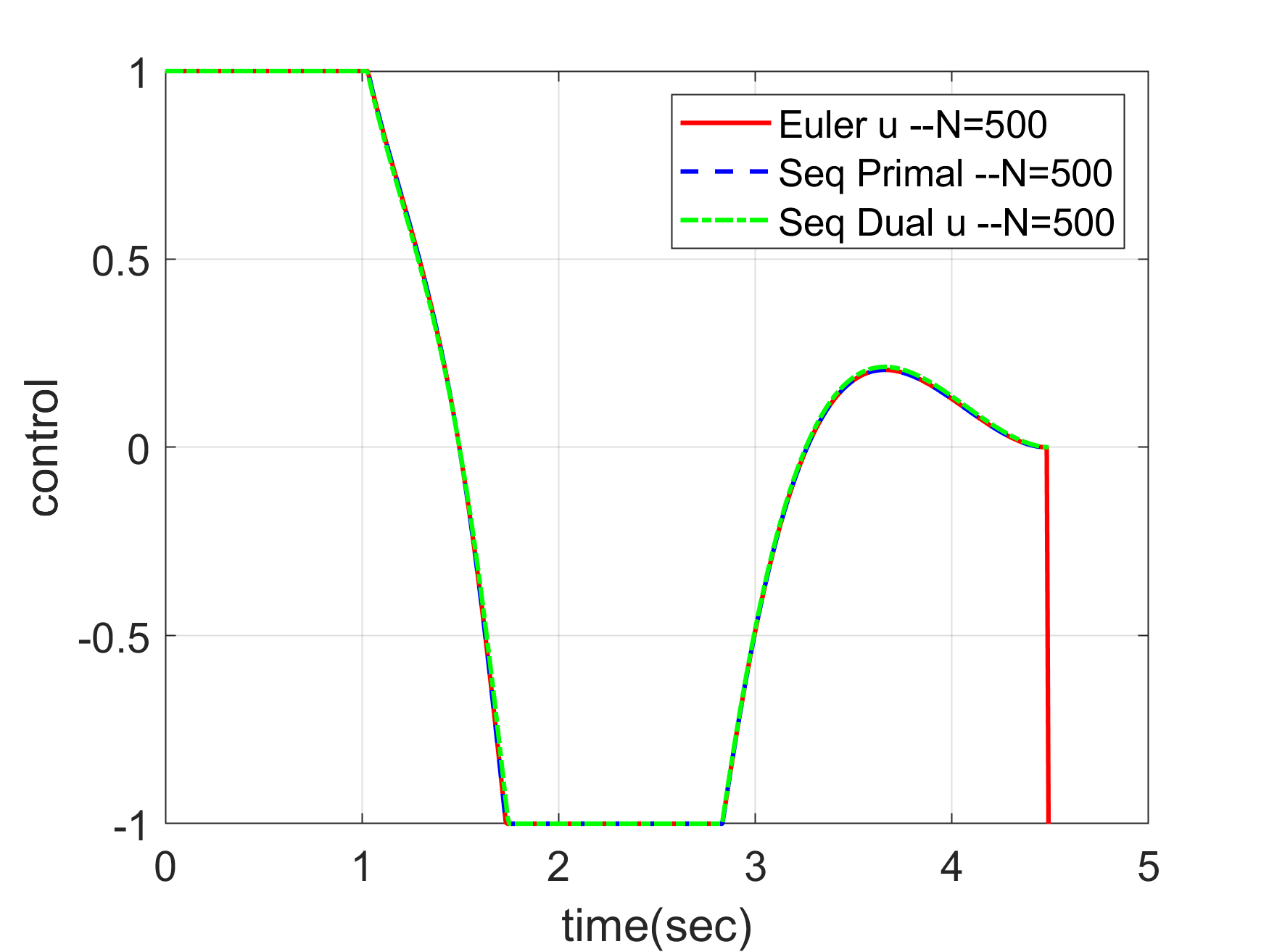}
		%	\caption{N=500时对偶方法的解与直接欧拉离散的解}
	\end{minipage}
	\begin{minipage}[t]{0.49\textwidth}
		\centering
		\includegraphics[width=6cm]{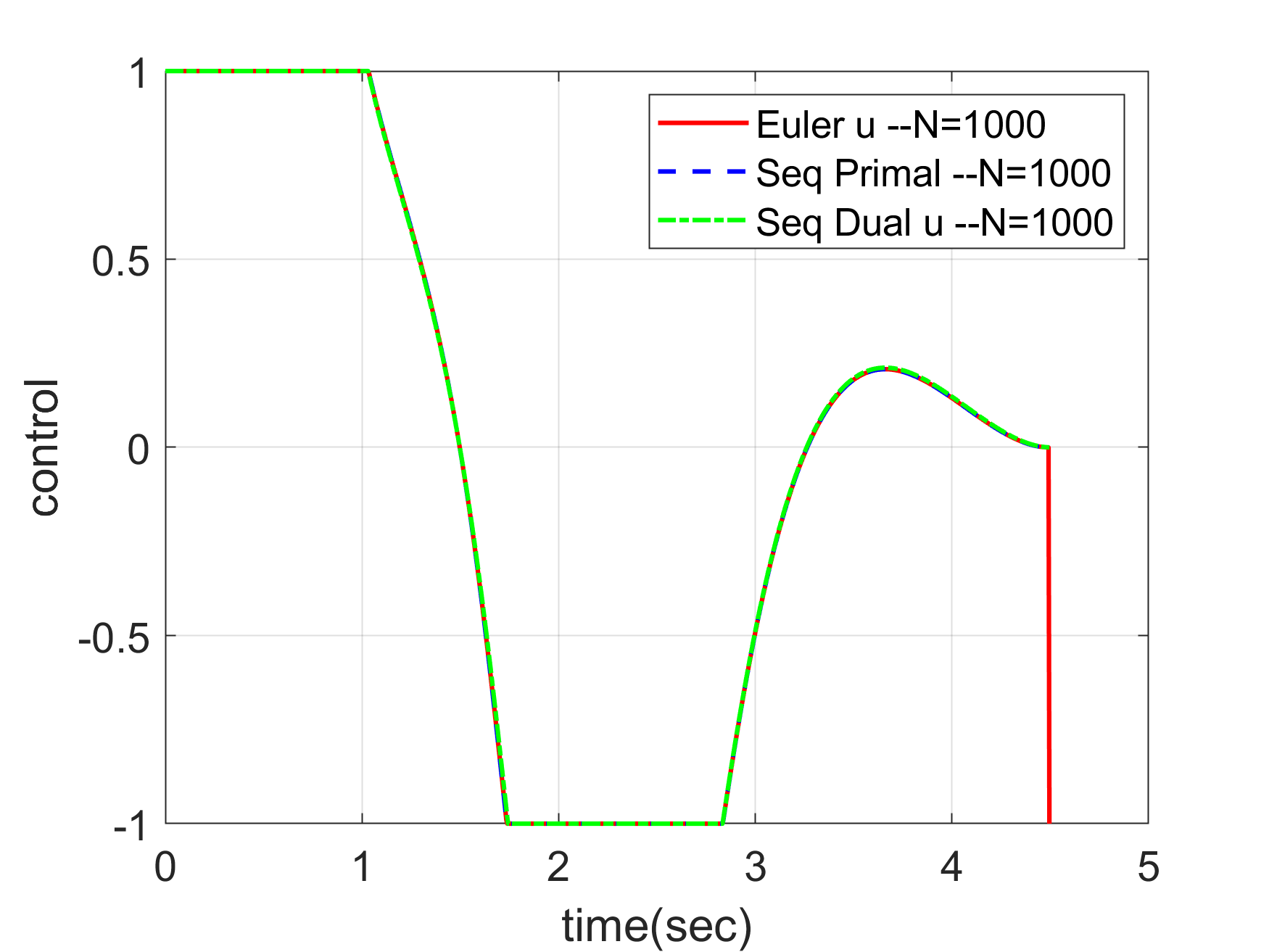}
		%	\caption{N=1000时对偶方法的解与直接欧拉离散的解}
	\end{minipage}
	\caption{Example 3 optimal control obtained with N=500 and N=1000}
\end{figure}

Example 4. Consider the following problem that terminal state is fixed which is adopted from \cite{quasi}.  
\begin{align*}
\min \quad &\frac{1}{2}\int_{0}^{5} x_1^{2}+x_2^2+u^2 dt\\
\text{subject to}\quad &\dot{x_1}=x_2\\
&\dot{x_2}=(1-x_1^2)\cdot x_2-x_1+u\\
&-0.75\leq u\leq 0.75 \\
& x(0)=(1\quad 0)^{\top},\quad x(5)=(-1\quad 0)^{\top}
\end{align*}
The nonlinear problem above is of the form (NP) with $E=I$ ($I$ means identity matrix) and $\bs e_f=(-1\quad 0)^{\top}$. We use Algorithm 2 to solve the problem, constant $\theta$ in merit function is set $\theta=100$ and the max iteration number for each line search is set $50$. The following table shows the computation results and our method behaves well camparing to other two method.
\begin{table}[H]
	\caption{Example 4 numerical performance }
	\begin{tabular}{| l ||c|c | c| c | c| c| c | c |}
			\hline 
		\multirow{2}{*}{N} & \multicolumn{3}{|c|}{Optimal Value} & \multicolumn{2}{|c|}{number of quasilinearization} & \multicolumn{3}{|c|}{IPOPT time[sec]} \\
		\cline{2-9}
		\multirow{2}{*}{~} &Euler & SeqPrimal & SeqDual & SeqPrimal & SeqDual & Euler & SeqPrimal & SeqDual \\
		\hline
		\hline
		50  & 2.1983 &2.1983 & 2.0978 &5 &4 &0.021& 0.136 &0.018 \\
		\hline
		100  & 2.1643 &2.1643 & 2.1145 &5 &4 &0.032& 0.163 &0.020 \\
		\hline
		200  & 2.1497 &2.1497 & 2.1249 &5 &4 &0.073& 0.206 &0.030 \\
		\hline
		500  & 2.1416&2.1416 & 2.1318 & 5 &4  &0.156&0.552 &0.054 \\
		\hline
		1000 & 2.1391&2.1391 & 2.1341 &5 &4 &0.269&0.724  &0.085  \\
		\hline
		2000  & 2.1378 &2.1378 & 2.1354 &5 &4 &0.495& 1.582 &0.119 \\
		\hline
		5000 & 2.1371& 2.1371 & 2.1361 &5 &4 &1.45&4.74 &0.325  \\
		\hline
		10000 & 2.1368 &2.1368 & 2.1363 & 5 & 4 &4.11&7.98 &1.648 \\
		\hline
	\end{tabular}
\end{table}

\begin{figure}[H]
	\centering
	\begin{minipage}[t]{0.5\textwidth}
		\centering
		\includegraphics[width=6cm]{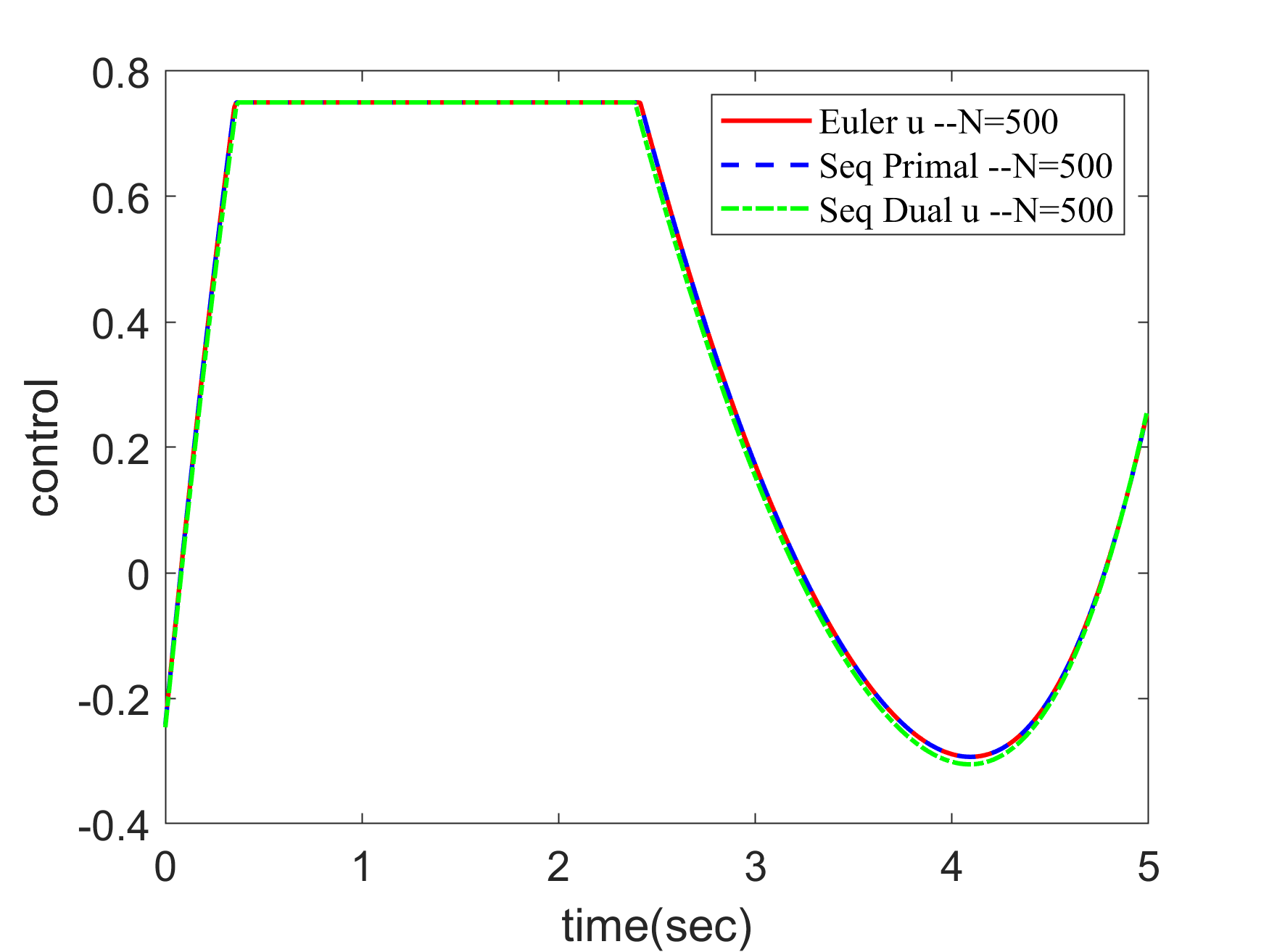}
		%	\caption{N=500时对偶方法的解与直接欧拉离散的解}
	\end{minipage}
	\begin{minipage}[t]{0.49\textwidth}
		\centering
		\includegraphics[width=6cm]{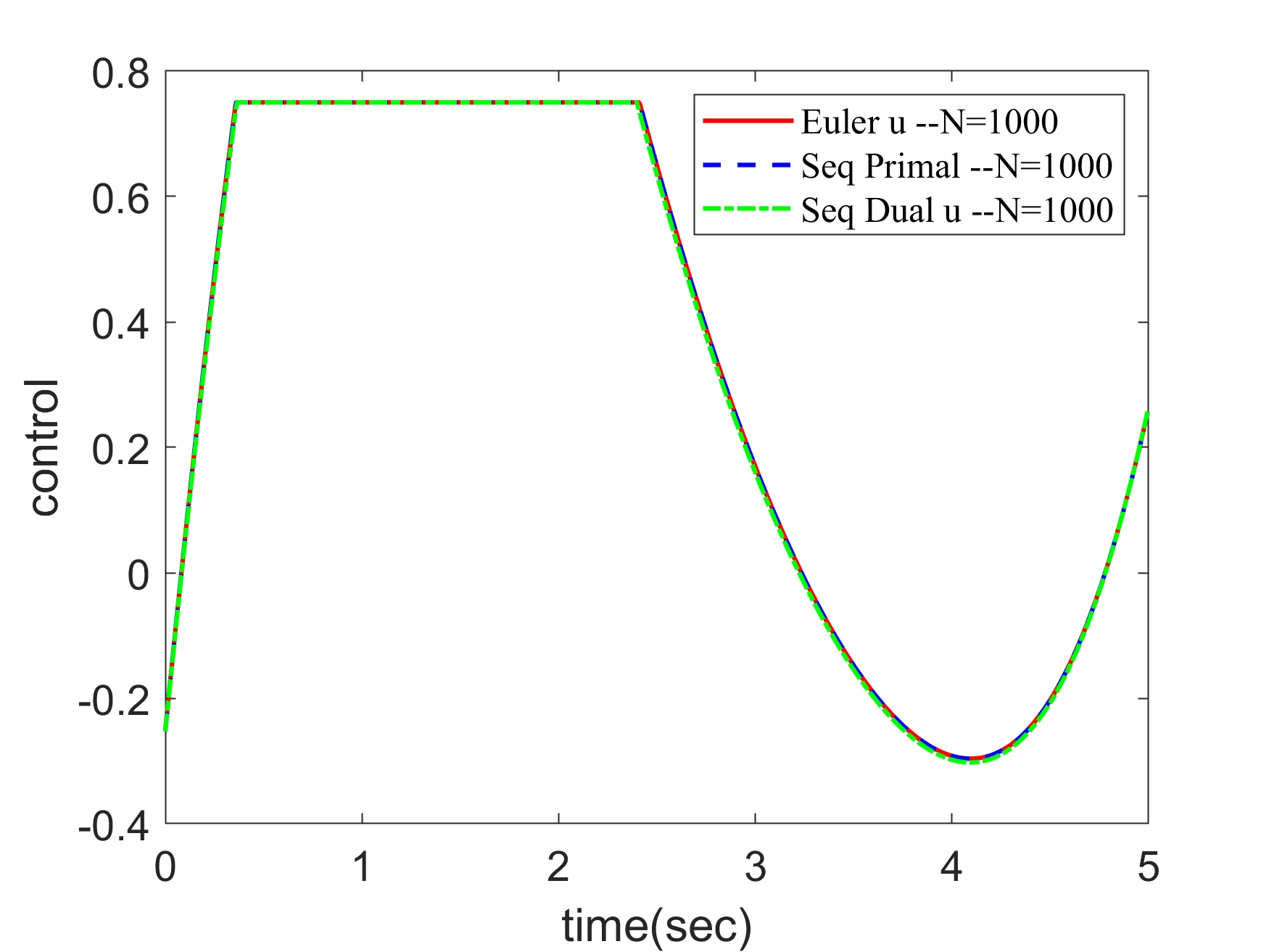}
		%	\caption{N=1000时对偶方法的解与直接欧拉离散的解}
	\end{minipage}
	\caption{Example 4 optimal control obtained with N=500 and N=1000}
\end{figure}

Example 5.The final example is adopted from \cite{quasi_bound}
\begin{align*}
\min \quad &\frac{1}{2}\int_{0}^{2.4} x_1^{2}+x_2^2+u^4+u^2 dt\\
\text{subject to}\quad &\dot{x_1}=x_2\\
&\dot{x_2}=(1-x_1^2)\cdot x_2-x_1+u\\
&-0.25\leq u\leq 1 \\
& x(0)=(1\quad 0)^{\top},\quad x(2.4)=(0\quad 0)^{\top} 
\end{align*}
We also notice that the nonlinear problem above is of the form (NP) with $E=I$ ($I$ means identity matrix) and $\bs e_f=(-1\quad 0)^{\top}$. We use Algorithm 2 to solve this problem, the constant $\theta$ in merit function and max iteration number for line search is chosen same as Example 4. The following table shows the computation results. Here the objective functional is strongly convex correspinding to state and control variables.
\begin{table}[H]
		\caption{Example 5 numerical performance }
\begin{tabular}{| l || c | c | c| c|  c | c|c |c|}
	\hline 
	\multirow{2}{*}{N} & \multicolumn{3}{|c|}{Optimal Value} & \multicolumn{2}{|c|}{number of quasilinearization} & \multicolumn{3}{|c|}{IPOPT time[sec]} \\
	\cline{2-9}
	\multirow{2}{*}{~} &Eluer & SeqPrimal & SeqDual  & SeqPrimal & SeqDual & Euler & SeqPrimal & SeqDual \\
	\hline
	\hline
	50 & 2.6005 &2.6004 &2.2743 &4 &5 &0.046 &0.137 &0.027\\
	\hline
	100 & 2.4969 &2.4969 &2.3398 &4 &5 &0.048 &0.192 &0.031 \\
	\hline
	200  & 2.4488&2.4488 & 2.3737 &5 & 5 &0.079&0.264 &0.045 \\
	\hline
	500  & 2.4244&2.4244 & 2.3947  &5 & 5  &0.153&0.511 &0.059  \\
	\hline
	1000 & 2.4167&2.4167 & 2.4019 &5 & 5&0.284&0.808 &0.102\\
	\hline
	2000 & 2.4130&2.4130 & 2.4056 &5 & 5&0.538&1.384 &0.115\\
	\hline
	5000 & 2.4107&2.4107 & 2.4078 &5 &5 &1.32&3.532 &0.327 \\
	\hline
	10000 & 2.4100&2.4100 & 2.4083 &5 &5 &2.82&7.06 &1.271 \\
	\hline
\end{tabular}
\end{table}
\begin{figure}[htbp]
	\centering
	\begin{minipage}[t]{0.49\textwidth}
		\centering
		\includegraphics[width=6cm]{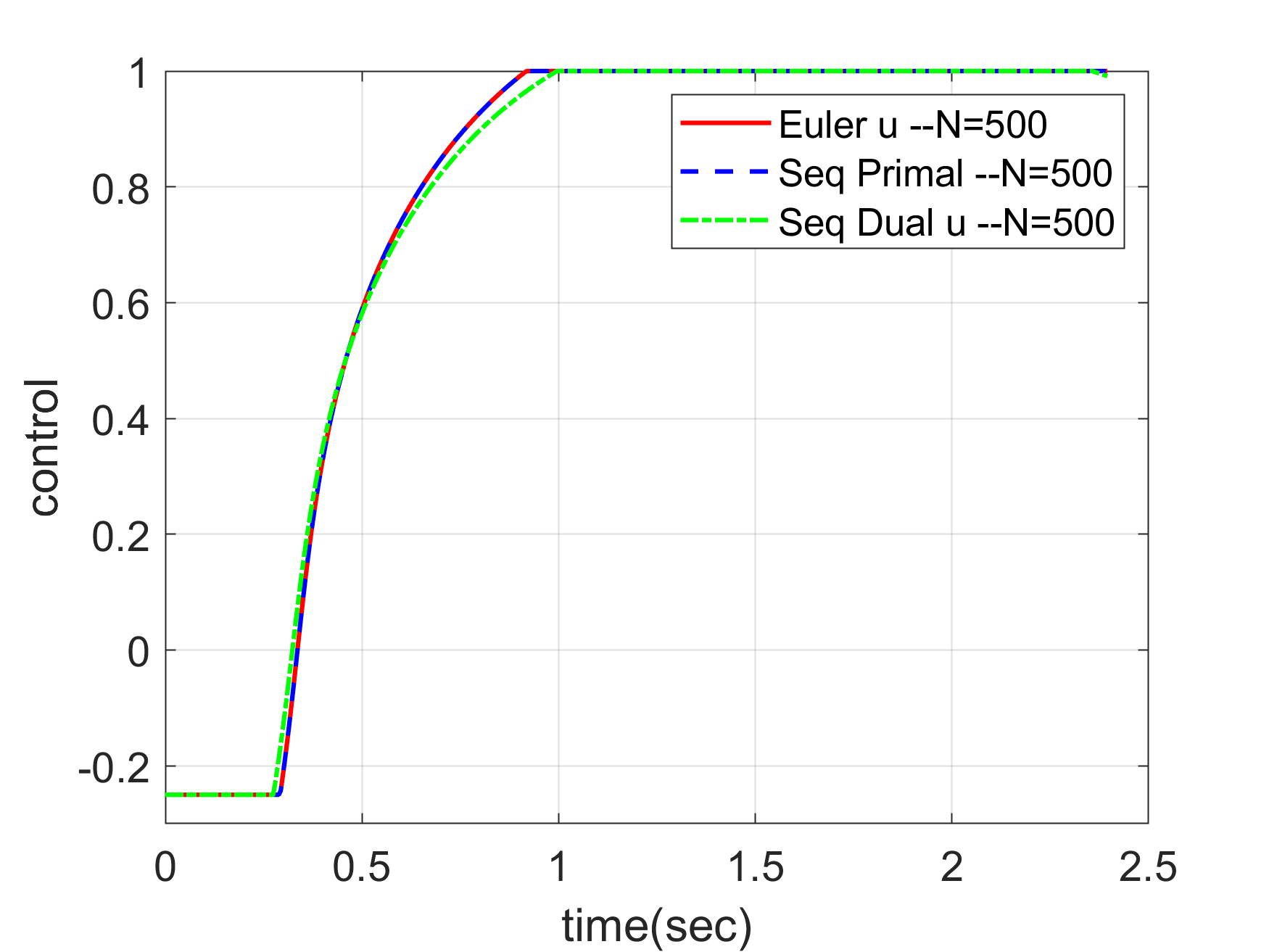}
		%	\caption{N=500时对偶方法的解与直接欧拉离散的解}
	\end{minipage}
	\begin{minipage}[t]{0.49\textwidth}
		\centering
		\includegraphics[width=6cm]{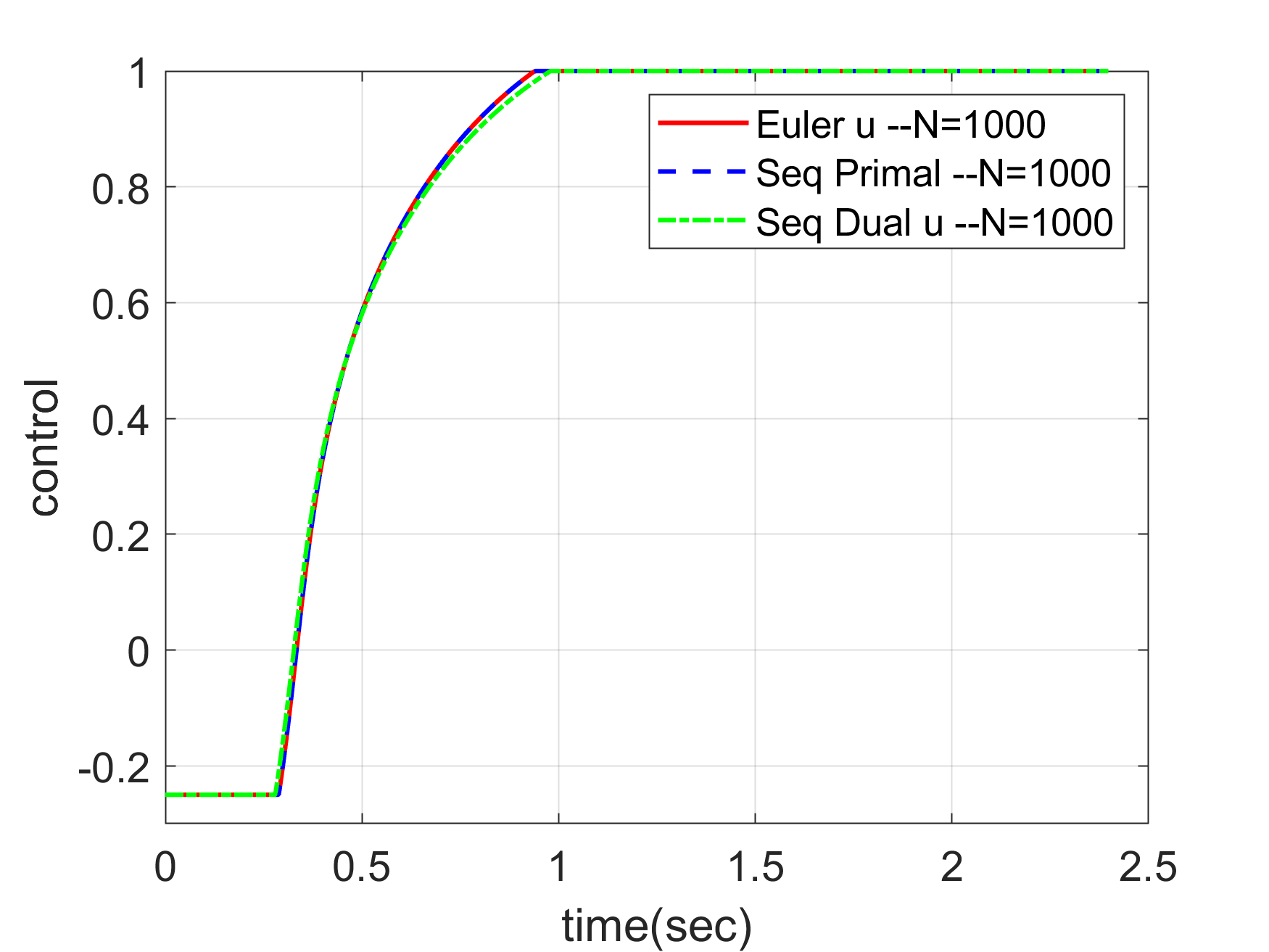}
		%	\caption{N=1000时对偶方法的解与直接欧拉离散的解}
	\end{minipage}
	\caption{Example 5 optimal control obtained with N=500 and N=1000}
\end{figure}

%\newpage

\section{Conclusion}
\qquad In this paper, we focused on a special class of nonlinear optimal control problem with control constraints and discussed how to solve it effectively. Firstly, we converted solving the original nonlinear problem into solving a sequence of linear-quadratic problem with control constraint whose solution is easier than the original nonlinear optimal control problem. Then we took advantage of Fenchel duality scheme to formulate the dual problem of the subproblem obtained by quasilinearization technique which could be described by a linear-quadratic problem. The strong duality and saddle point properties were analyzed in this process, thus it showed solving dual problem is effective. Our iterative scheme replaced solving each subproblem by solving its dual problem. Besides the convergence result of the algorithm designed for solving this class problem was also analyzed. Preliminary numerical results were reported to verify the theoretical assertions including the convergence of algorithm, strong duality for the control-constrained linear-quadratic optimal control problem as well as the convergence of Euler discretization. The numerical results also demonstrated that the proposed method were computationally efficient and easy to implement without sacrificing the accuracy of the solution. Future work should consider more general problems with convex objective functional for example the case when $g=0$. Besides more general constraints such as pointwise state constraints and so on should also be included in our furture work. 

\clearpage
\bibliography{refee}

\end{document}